\documentclass[smallextended]{svjour3}       
\usepackage{graphicx}
\usepackage{amsmath}
\usepackage{amssymb}
\usepackage{color}
\usepackage[letterpaper, margin=1.5in]{geometry}
\newcommand{\cci}{C} 
\newcommand{\set}[1]{\left\{#1\right\}}                     
    
\begin{document}

\title{Two-sided linear chance constraints and extensions
}


\author{Miles Lubin  \and
        Daniel Bienstock \and
        Juan Pablo Vielma
}


\institute{M. Lubin \and J. P. Vielma \at
              Massachusetts Institute of Technology, Cambridge, MA, USA \\
              \email{mlubin@mit.edu}           
           \and
           D. Bienstock \at
           Columbia University, New York, NY, USA
}

\date{February 2016}

\maketitle

\begin{abstract}
We examine the convexity and tractability of the two-sided linear chance constraint model under Gaussian uncertainty. We show that these constraints can be applied directly to model a larger class of nonlinear chance constraints as well as provide a reasonable approximation for a challenging class of quadratic chance constraints of direct interest for applications in power systems. 
With a view towards practical computations, we
develop a second-order cone outer approximation of the two-sided chance constraint with provably
small approximation error.
\end{abstract}

\section{Introduction}

Chance constraints (or probabilistic constraints) were among the first extensions proposed to linear programming as a natural formulation for treating constraints where some of the coefficients are uncertain at the time of optimization~\cite{CharnesCooper}. In the chance constraint model, we suppose that the uncertain values follow a known distribution and enforce that the constraint holds with high probability as a function of the decision variables.

Nemirovski and Shapiro~\cite{NS2007} observe that, in general, convexity and tractability results in chance constraints are a rare combination. When the corresponding deterministic constraint is convex, the chance constraint may be nonconvex. And even for those chance constraints which are in fact convex, the authors \cite{NS2007} cite examples where such constraints remain computationally intractable because it is NP-Hard to test if the constraint is satisfied. For linear chance constraints of the form
\begin{equation}\label{eq:onesidechance}
    \mathbb{P}(x^T\xi \le b) \ge 1-\epsilon,
\end{equation}
where $x \in \mathbb{R}^n$ and $b \in \mathbb{R}$ are decision variables, the constraint is known to be convex (that is, the set $\{ (x,b) : \mathbb{P}(x^T\xi \le b) \ge 1-\epsilon\}$ is convex)  and computationally tractable when $\xi$ has an \textit{elliptical log-concave} distribution~\cite{Lagoa05}, examples of which include the multivariate Gaussian distribution and few others. The computational challenges presented by chance constraints have motivated approximation schemes~\cite{NS2007} and alternative formulations such as robust optimization~\cite{BenTalNemirovskiRobust2000}.

Even more challenging than linear chance constraints, \textit{joint chance constraints} require that a set of linear constraints hold jointly with high probability. Pr\'{e}kopa~\cite{PrekopaBook} reviews many of the standard results. In particular, he proves convexity of the constraint $\mathbb{P}(x \ge \xi) \ge 1-\epsilon$ with respect to $x\in \mathbb{R}^n$ when $\xi$ follows a multivariate continuous log-concave distribution and of the constraint $\mathbb{P}(Tx \ge 0) \ge 1-\epsilon$ when some elements of the matrix $T$ are random with a joint Gaussian distribution and have a specialized covariance structure between the rows of $T$ (further generalized by~\cite{Copulas}). Van Ackooij et al.~\cite{VanAckooij10} consider \textit{rectangular} chance constraints of the form $\mathbb{P}(a \le \xi \le b) \ge 1-\epsilon$ with respect to vectors $a$ and $b$ where $\xi$ follows a multivariate Gaussian distribution. Their model does not allow for products between random variables and decision variables.

The basic model we consider in this work, which is a special case of a joint chance constraint, is the two-sided chance constraint
\begin{equation}\label{eq:twosideintro}
    \mathbb{P}(a \le x^T\xi \le b) \ge 1-\epsilon,
\end{equation}
where $a \in \mathbb{R}, b \in \mathbb{R}$, and $x \in \mathbb{R}^n$ are decision variables, and $\xi$ is jointly Gaussian with known mean and covariance. In Section~\ref{sec:cvx2side}, we prove that this constraint is in fact convex in $a$, $b$ and $x$ given $\epsilon \le \frac{1}{2}$. The proof, which we believe is the first, follows from a geometrical insight combined with standard tools for chance constraints such as log-concavity. The major methodological contributions of this work lie in the subsequent generalizations of the model and in our analysis of the computational tractability of the chance constraint. In Section~\ref{sec:exact_extensions} we show that a number of seemingly more complex and nonlinear constraints can be formulated by using the two-sided constraint~\eqref{eq:twosideintro}. In Section~\ref{sec:tractability}, we demonstrate computational tractability of these constraints under a modern mathematical optimization lens. In addition to an exact derivative-based nonlinear formulation, we develop an approximate second-order cone (SOC) formulation for~\eqref{eq:twosideintro} with provable approximation quality. This SOC formulation permits one to incorporate such constraints into large-scale models solvable by state-of-the-art commercial and open-source software.

Using~\eqref{eq:twosideintro} as a primitive, we develop an approximation for the more challenging chance constraint
\begin{equation}\label{eq:quadintro}
\mathbb{P}((a^T\xi + b)^2 + (c^T\xi + d)^2 \le k) \ge 1-\epsilon,
\end{equation}
where $a,c \in \mathbb{R}^n$, $b,d,k \in \mathbb{R}$ are decision variables, and $\xi$ is jointly Gaussian with known mean and covariance. This constraint is motivated by applications in power systems which we discuss in Section~\ref{sec:motivation}. In Section~\ref{sec:quadapprox}, we study the constraint~\eqref{eq:quadintro} in detail and compare a number of approximation schemes, ultimately demonstrating that our approximation based on two-sided constraints is reasonable and of practical interest for its tractability.

\section{Motivation}\label{sec:motivation}

The basic question which motivates this work is the  short-term planning problem, known as \textit{optimal power flow} (OPF), which is solved as part of the real-time operation of the power grid to determine the minimum-cost production levels of controllable generators subject to reliably delivering electricity to customers across a large geographical area~\cite{OPFreview,BergenBook}. Conceptually, OPF is similar to a network flow problem with the additional complication that power flows according to the nonlinear Kirchhoff laws. On top of the nonlinear power flow laws, we aim to consider the uncertainty in production levels of renewable energy sources such as wind and solar photovoltaic.

In its traditional, deterministic form, OPF seeks to minimize total production costs
\begin{equation}
\operatorname*{minimize}_{p,\theta,f} \sum_{i \in \mathcal{G}} c_{i}p_i\label{eq:Det_OPF}
\end{equation}
\text{subject to the constraints} 
\begin{gather}
\sum_{n : \{b,n\} \in \mathcal{L}} f_{bn} - \sum_{m : \{m,b\} \in \mathcal{L}} f_{mb} = \sum_{i \in G_b} p_i + w_b - d_b, \quad \forall b \in \mathcal{B},\label{eq:balance} \\
\label{eq:gencapacity} p_{i}^{min} \leq p_i \leq p_i^{max}, \quad \forall i \in \mathcal{G},\\
\label{eq:flowdef} f_{mn} = \beta_{mn}(\theta_m - \theta_n), \quad \forall \{m,n\} \in \mathcal{L}, \\
\label{eq:Det_OPF_end} -f_{mn}^{max} \leq f _{mn}  \leq f_{mn}^{max}, \quad \forall \{m,n\} \in \mathcal{L}, 
\end{gather}
where $\mathcal{B}$ is the set of nodes (buses) in the grid, $\mathcal{G}$ is the set of generators, $G_b$ is the set of generators located at node $b$, and $\mathcal{L}$ is the set of edges (transmission lines). Decision variables $p_i$ denote the production levels of generator $i$, and the variables $f_{mn}$ denote the flow from node $m$ to node $n$. The value $d_b$ is the demand at each node (assumed to be known), and the value $w_b$ is the forecast production level from renewable energy sources (again assumed to be known). Constraint~\eqref{eq:balance} is the familiar flow balance constraint which balances supply with demand at each node. Constraints~\eqref{eq:gencapacity} and \eqref{eq:Det_OPF_end} enforce the capacities of the generators and transmission lines, respectively. The constraint~\eqref{eq:flowdef} links the flows to the bus angles $\theta$ and arises from the standard ``DC'' linearization of the nonlinear power flow laws; hence, this formulation is often called DCOPF. The formulation as stated above is efficiently solvable by linear programming on large-scale systems with tens of thousands of nodes within real-time operational constraints.

Our motivation is to address two major deficiencies in the standard DCOPF model. The first major deficiency is the deterministic nature of the model. In particular, the amount of power generated by renewable energy sources such as wind is highly variable and must be accounted for in short-term planning.

The line of work by~\cite{ccopf-sirev,JuMPChanceCaseStudy} addresses this deficiency by introducing chance constraints. More specifically, Bienstock et al.~\cite{ccopf-sirev} propose to model the deviations from the forecast wind production levels as zero-mean Gaussian random variables $\boldsymbol\omega_b$, combined with a proportional response policy for the generators. Letting $\Omega$ be the total, real-time deviation from the forecast (a positive value if there is more renewable generation than expected), each generator has a proportional response coefficient $\alpha_i$ and adjusts its real-time production to match $p_i - \alpha_i\Omega$. If $\sum_i \alpha_i = 1$, then this response policy guarantees balance of supply and demand, although it does not guarantee that output capacities or transmission capacities are always satisfied. Both $p_i$ and $\alpha_i$ are decision varibles. Transmission capacities, in practice, are soft constraints, and hence~\cite{ccopf-sirev} propose to enforce them as chance constraints
\begin{equation}\label{eq:chanceabsflow}
\mathbb{P}(|\boldsymbol f_{mn}| \le f_{mn}^{max}) \ge 1-\epsilon,
\end{equation}
where $\boldsymbol f_{mn}$ is the random flow driven by the deviations $\boldsymbol\omega_b$. Bienstock et al.~\cite{ccopf-sirev} then approximate~\eqref{eq:chanceabsflow} by splitting it into two constraints
\begin{equation}
\mathbb{P}(\boldsymbol f_{mn} \le f_{mn}^{max}) \ge 1-\epsilon \text{ and } \mathbb{P}(\boldsymbol f_{mn} \ge -f_{mn}^{max}) \ge 1-\epsilon,
\end{equation}
both of which can be expressed as simple linear Gaussian chance constraints~\eqref{eq:onesidechance}. The assumption that deviations from the forecast follow a Gaussian distribution is made for tractability. This assumption can be further refined, with practical gain, without loss of tractability by introducing uncertainty sets on the parameters of the Gaussian distribution~\cite{JuMPChanceCaseStudy}.

The second major deficiency in the standard DCOPF model is the crude approximation it provides of the true, nonlinear, nonconvex power flow laws. In particular, the linearized model assumes constant voltage and therefore neglects so-called \textit{reactive} power flow, which is the imaginary component of complex-valued power flow. The real component is referred to as \textit{active} power. Although we cannot directly treat the nonconvex case, we propose to consider more accurate linearizations which account for reactive power, such as those which arise from linearizing around a current operating solution~\cite{bolognani2015fast}. When extending the model of~\cite{ccopf-sirev} to account for reactive power, we obtain chance constraints of the form
\begin{equation}\label{eq:quadpower}
\mathbb{P}((\boldsymbol f_{mn}^{active})^2 +(\boldsymbol f_{mn}^{reactive})^2  \le (f_{mn}^{max})^2) \ge 1-\epsilon,
\end{equation}
because transmission capacities are limited by the magnitude of the complex-valued power flow across a line.

Our first attempt at studying the constraint~\eqref{eq:quadpower} led us to study the simpler two-sided form~\eqref{eq:twosideintro}. These results, in turn, provided us with a means to approximate~\eqref{eq:quadpower}, as we discuss in Section~\ref{sec:quadapprox}. The approximation we derive here has already yielded a practical implementation in the JuMPChance modeling package~\cite{JuMPChance} which is being used to study the value of the model we propose in ongoing work~\cite{working}.

\section{Convexity of two-sided Gaussian linear chance constraints}\label{sec:cvx2side}

The main result in this section is the convexity of the two-sided chance constraint~\eqref{eq:twosideintro}.

Let $\varphi(x) = \frac{1}{\sqrt{2\pi}}e^{-\frac{x^2}{2}}$ be the standard Gaussian density and $\Phi(x) = \int_{-\infty}^x \varphi(t)\,dt$ the Gaussian integral.

\begin{definition}
    Let $\xi \sim N(0,1)$ be a standard Gaussian random variable. Let $\epsilon \in (0,1)$. We define the set $S_\epsilon := \{ (x,y) \in \mathbb{R}^2 : \mathbb{P}(x \le \xi \le y) \ge 1-\epsilon \}.$
\end{definition}
Note that $S_\epsilon$ has two equivalent representations as $\{ (x,y) : \int_x^y \varphi(t)\, dt \ge 1-\epsilon \}$ and $\{ (x,y) : \Phi(y) - \Phi(x) \ge 1-\epsilon \}$.

We will proceed to prove that $S_\epsilon$ is convex, but first we define \textit{log-concavity} and recall some basic properties. See Boyd~\cite{BoydBook} for further discussion and proofs of these properties.
\begin{definition}
A non-negative function $f : \mathbb{R}^n \to \mathbb{R}$ is log-concave if $\forall\, x,y \in \operatorname{dom} f \text{ and } \lambda \in (0,1)$
\[
    f(\lambda x + (1-\lambda)y) \geq f(x)^\lambda f(y)^{1-\lambda}.
\]
\end{definition}

For strictly positive functions $f$, this definition is equivalent to the condition that $\log f$ is concave. It is easy to verify, therefore, that the Gaussian density $\varphi$ is log-concave. Lemma~\ref{lem:logconcave} recalls basic properties of log-concave functions.

\begin{lemma}\label{lem:logconcave} The following properties hold for log-concave functions:
\begin{itemize}
    \item If $f : \mathbb{R}^n \to \mathbb{R}$ and $g : \mathbb{R}^n \to \mathbb{R}$ are log-concave, then the product $h(x) = f(x)g(x)$ is log-concave.
    \item If $f : \mathbb{R}^n \to \mathbb{R}$ is the indicator function of a convex set, then $f$ is log-concave.
    \item If $f : \mathbb{R}^n \times \mathbb{R}^m \to \mathbb{R}$ is log concave, then $g(x) = \int f(x,y)\,dy$ is log-concave on $\mathbb{R}^n$.
\end{itemize}
\begin{proof}
    See Boyd~\cite{BoydBook}.
\end{proof}
\end{lemma}

With these basic properties, we can proceed to prove the following lemma.

\begin{lemma}\label{lem:seps}
The set $S_\epsilon$ is convex.

\begin{proof}
    Let $I(s,r) = 1$ if $s \le r$ and zero otherwise. That is, $I$ is the indicator function for the convex set $\{(a,b) : a \le b\}$. Therefore the function $g(t,x,y) = \varphi(t)I(t,y)I(x,t)$ is log concave, because it is a product of log concave functions. Then for $y \geq x, f(x,y) = \int_x^y \varphi(t)\,dt = \int \varphi(t)I(t,y)I(x,t)\,dt$ is log concave, because it is the marginal of a log concave function. Hence $S_\epsilon$ is convex because it is an upper level set of a log-concave function.
\end{proof}
\end{lemma}

Convexity of $S_\epsilon$ proves convexity of the very simple chance constraint $\mathbb{P}(x \le \xi \le y) \ge 1-\epsilon$ for all $\epsilon \in (0,1)$ with respect to $(x,y) \in \mathbb{R}^2$. Note that this convexity result is a special case of the rectangular constraints considered by~\cite{VanAckooij10}. In order to account for products between the decision variables and the random variables, we require the following additional developments.

\begin{definition}
Let $\bar S_\epsilon = \operatorname{cl} \{ (x,y,z) : (x/z,y/z) \in S_\epsilon, z > 0 \}$ be the conic hull of $S_\epsilon$ (where $\operatorname{cl}$ is the closure operator).
\end{definition}

By standard results~\cite{HiriartLemarechal93book2}, $\bar S_\epsilon$ is convex. The following lemma, in which we prove monotonicity properties of the set $\bar S_\epsilon$, is key to our main result.

\begin{lemma}
\label{lem:monotonicity}
Let $\epsilon \in (0,\frac{1}{2}]$. Then $(x,y,z) \in \bar S_\epsilon$ iff $z \ge 0$ and $\exists\, x' \ge x, y' \le y,$ and $z' \ge z$ such that $(x',y',z') \in \bar S_\epsilon$.

\begin{proof}
    Suppose we are given $(x',y',z') \in \bar S_\epsilon$ and $(x,y,z)$ with $x \le x', y \ge y',$ and $0 < z \le z'$. We will show that $(x,y,z) \in \bar S_\epsilon$.   By symmetry of the Gaussian density and $\epsilon \le \frac{1}{2}$, $(x',y',z') \in \bar S_\epsilon$ implies $x' < 0$ and $y' > 0$, so $x/z \le x'/z \le x'/z'$ and $y/z \ge y'/z \ge y'/z'$. By increasing the upper limit of integration or decreasing the lower limit of integration, we can only increase the value of the integral, so
\begin{equation}
\int_{x/z}^{y/z} \phi(t)\, dt \ge \int_{x'/z'}^{y'/z'} \phi(t)\, dt \ge 1-\epsilon.
\end{equation}
For the case of $z=0$, take a sequence of decreasing iterates $z_1 = z', z_2, z_3, \ldots$ with $z_i \to 0$. For each $i$, the above argument shows $(x,y,z_i) \in \bar S_\epsilon$, which implies $(x,y,0) \in \bar S_\epsilon$ since $\bar S_\epsilon$ is a closed set.
\end{proof}
\end{lemma}

With these properties, we now prove the main result of this section.

\begin{theorem}\label{thm:chance}
Let $\xi$ be a vector of $n$ i.i.d.\ standard Gaussian random variables, $0 < \epsilon \le \frac{1}{2}$ and 
\[\cci := \set{(a,b,x) \in \mathbb{R}\times\mathbb{R}\times\mathbb{R}^n\,:\, \mathbb{P}(a \le x^T\xi \le b) \geq 1 - \epsilon }. \]
Then $\cci$ is a projection of the convex set 
\[ \set{(a,b,x,t) \in \mathbb{R}\times\mathbb{R}\times\mathbb{R}^n\times \mathbb{R}\,:\,  ||x||_2 \leq t,\quad (a,b,t) \in \bar S_\epsilon }, \] 
and hence $\cci$ is convex. 
\begin{proof}
\begin{equation}
\mathbb{P}(a \leq x^T\xi \leq b) \geq 1-\epsilon
\end{equation}
iff
\begin{equation}\label{eq:probrange}
\mathbb{P}\left(\frac{a}{||x||_2} \leq \frac{x^T\xi}{||x||_2} \leq \frac{b}{||x||_2}\right) \geq 1-\epsilon
\end{equation}
iff
\begin{equation}\label{eq:transform_to_s_eps}
(a,b,||x||_2) \in \bar S_\epsilon
\end{equation}
iff (by Lemma~\ref{lem:monotonicity})
\begin{equation}
\exists\, t \ge ||x||_2 \text{ such that } (a,b,t) \in \bar S_\epsilon.
\end{equation}

Where the equivalence between \eqref{eq:probrange} and \eqref{eq:transform_to_s_eps} holds because $\frac{x^T\xi}{||x||_2}$ is a standard Gaussian random variable. The above proof assumes $x \neq 0$. For the case of $x = 0$,
\begin{equation}
\mathbb{P}(a \le x^T\xi \le b) \geq 1 - \epsilon
\end{equation}
iff
\begin{equation}
a \le 0 \le b
\end{equation}
iff
\begin{equation}
(a,b,0) \in \bar S_\epsilon.
\end{equation}
The justification for the final equivalence is as follows. If the strict inequality $a < 0 < b$ holds, then $\lim_{z \to 0+} \int_{a/z}^{b/z} \varphi(t)\, dt = 1$, so membership holds in $\bar S_\epsilon$. If $a=0$, $b=0$, or both, then we can construct a sequence of points $(a_i,b_i,0) \to (a,b,0)$ with each $(a_i,b_i,0) \in \bar S_\epsilon$, so the statement holds because $\bar S_\epsilon$ is closed.
\end{proof}
\end{theorem}

More generally,

\begin{lemma}\label{lem:generalcc}
Let $\xi \sim N(\mu,\Sigma)$ be a jointly distributed Gaussian random vector with mean $\mu$ and positive definite covariance matrix $\Sigma$ and $0 < \epsilon \le \frac{1}{2}$, and let 
\[\cci_{\mu,\Sigma} := \set{(a,b,x) \in \mathbb{R}\times\mathbb{R}\times\mathbb{R}^n\,:\, \mathbb{P}(a \le x^T\xi \le b) \geq 1 - \epsilon }. \]
Then $\cci_{\mu,\Sigma}$ convex.

\begin{proof}
Let $LL^T = \Sigma$ be the Cholesky decomposition of the covariance matrix $\Sigma$. Then $\xi = L\zeta + \mu$ where $\zeta$ is a vector of i.i.d.\ standard Gaussian random variables.
The point $(a,b,x)$ satisfies
\begin{equation}
\mathbb{P}(a \leq x^T\xi \leq b) \geq 1-\epsilon
\end{equation}
iff
\begin{equation}
\mathbb{P}(a \leq x^T(L\zeta + \mu) \leq b) \geq 1-\epsilon
\end{equation}
iff
\begin{equation}
(a-\mu^Tx,b-\mu^Tx,L^Tx) \in C.
\end{equation}
That is, the set $\cci_{\mu,\Sigma}$ is an affine transformation of the convex set $C$ representing the i.i.d.\ case, and hence $\cci_{\mu,\Sigma}$ is convex.

\end{proof}
\end{lemma}

\section{Exact extensions of two-sided constraints}\label{sec:exact_extensions}

In this section, we generalize the basic result in Section~\ref{sec:cvx2side} to a number of cases in which a seemingly more complex chance constraint can be represented exactly by using two-sided chance constraints.

\subsection{Nonlinear chance constraints}

The simplest nonlinear constraint we consider, which will be used in formulating the approximation of the quadratic chance constraint in Section~\ref{sec:quadapprox}, is the absolute value constraint.

\begin{lemma}\label{eq:absconvex}
Let $\xi \sim N(\mu,\Sigma)$ be a jointly distributed Gaussian random vector with mean $\mu$ and positive definite covariance matrix $\Sigma$ and $0 < \epsilon \le \frac{1}{2}$.
Then the set
\begin{equation}
    \{ (a,b,x) \in \mathbb{R}\times\mathbb{R}\times\mathbb{R}^n: \mathbb{P}(|x^T\xi + a| \leq b) \geq 1 - \epsilon \}
\end{equation} 

is convex.

\begin{proof}
$\mathbb{P}(|x^T\xi + a| \leq b) \geq 1 - \epsilon$ \text{ iff } $\mathbb{P}(-b - a \le x^T\xi \le b - a) \geq 1 - \epsilon$.
\end{proof}
\end{lemma}

The above lemma is a special case of the following significantly more general theorem:

\begin{theorem}\label{thm:nonlinearchance}
Let $f: \mathbb{R} \to \mathbb{R}$ be a convex function which attains its minimum at $x = c$, let $g : \mathbb{R}^m \to \mathbb{R}$ be an arbitrary convex function, and let $\xi$ be a standard Gaussian random vector (without loss of generality, we can assume independence and zero mean). Let $\epsilon \le \frac{1}{2}$. Then the set
\begin{equation}\label{eq:nonlinconstr}
D := \left\{ (x,z,b) \in \mathbb{R}^n\times\mathbb{R}^m\times \mathbb R  : \mathbb{P}(f(x^T\xi + b) + g(z) \leq 0) \geq 1-\epsilon \right\}
\end{equation}
is a projection of the convex set
\begin{align}
    \{ (x,z,b,k,x',y',t) \in \mathbb{R}^{n+m+5}: &t \ge ||x||_2, k \le -g(z), 
    x' \ge l(k) - b - c, \\&y' \le u(k) - b - c,
    (x',y',t) \in \bar S_\epsilon \}
\end{align}
where $l$ and $u$, are explicitly computable convex and concave functions, respectively, which we define below depending on $f$. And hence, $D$ is convex.
\begin{proof}
Let $l(k)$ and $u(k)$ be functions such that $f(x-c) \leq k$ iff $x \in [l(k),u(k)]$. We can obtain $l$ and $u$ by shifting the graph of $f$ so that the minimum is at zero and then reflecting the graph along $y = x$, and since $f(\cdot-c)$ is decreasing up to zero and increasing after zero, we have in particular that $u(k)$ is concave and increasing and $l(k)$ is convex and decreasing.
Then
\begin{equation}
\mathbb{P}(f(x^T\xi + b) + g(z) \leq 0) \geq 1-\epsilon
\end{equation}
iff
\begin{equation}
\mathbb{P}(l(-g(z)) \leq x^T\xi + b + c \leq u(-g(z))) \geq 1-\epsilon
\end{equation}
iff
\begin{equation}
\mathbb{P}(l(-g(z)) - b - c \leq x^T\xi \leq u(-g(z)) - b - c) \geq 1-\epsilon
\end{equation}
iff (by Theorem~\ref{thm:chance})
\begin{equation}\label{eq:nonlinear_extended}
\exists\, t \ge ||x||_2 \text{ and } x' \ge l(-g(z)) - b - c \text{ and }  y' \le u(-g(z)) - b - c \text{ such that } (x',y',t) \in \bar S_\epsilon.
\end{equation}
Finally,
\begin{equation}
x' \ge l(-g(z)) - b - c \text{ and }  y' \le u(-g(z)) - b - c
\end{equation}
iff (by $l$ decreasing and $u$ increasing)
\begin{equation}\label{eq:nonlinear_extended2}
\exists\, k \le -g(z) \text{ such that } x' \ge l(k) - b - c \text{ and }  y' \le u(k) - b - c.
\end{equation}
Since $l(k)$ is convex and $u(k)$ is concave, conditions~\eqref{eq:nonlinear_extended} and \eqref{eq:nonlinear_extended2} give a convex formulation, in an extended set of variables, for the chance constraint~\eqref{eq:nonlinconstr}.
\end{proof}
\end{theorem}

Theorem~\ref{thm:nonlinearchance} is sufficiently general to shed light on the quadratic chance constraint~\eqref{eq:quadpower} which motivated our original work. If one of the terms in the chance constraint is deterministic, then the constraint is indeed convex, as the following lemma shows. This simpler form of the quadratic constraint itself can be useful for the motivating application in power systems, if, for example, the reactive power flow across a transmission line is not subject to randomness.

\begin{lemma}
Let $\xi \sim N(\mu,\Sigma)$ be a jointly distributed Gaussian random vector with mean $\mu$ and positive definite covariance matrix $\Sigma$ and $0 < \epsilon \le \frac{1}{2}$.
Then the set
\begin{equation}
    \{ (x,b,k,z) \in \mathbb{R}^n\times\mathbb{R}\times\mathbb{R}\times\mathbb{R}: \mathbb{P}((x^T\xi + b)^2 + z^2 \leq k) \geq 1 - \epsilon \}
\end{equation} 

is convex.
\begin{proof}
Set $f(y) = y^2$, $g(z,k) = z^2 - k$ and apply Theorem~\ref{thm:nonlinearchance}.
\end{proof}
\end{lemma}

A similar proof technique as used in Theorem~\ref{thm:nonlinearchance} can also be applied in other cases. In the following lemma, we demonstrate convexity of the quadratic chance constraint in another special case when the random variable is univariate.

\begin{lemma}
Let $\xi$ be a scalar standard Gaussian random variable and let $\epsilon \in (0, 1)$. Then the set
\[
\left\{ (b,d,k) \in \mathbb{R}^3 : \mathbb{P}((\xi+b)^2 +(\xi+d)^2 \leq k) \geq 1-\epsilon \right\}
\]
is convex.
\begin{proof}
By applying the quadratic formula, we see that $(\xi+b)^2 +(\xi+d)^2 \leq k$ iff $\xi \in [l(b,d,k),u(b,d,k)]$ where
\[
l(b,d,k) = \frac{1}{2}\left(-(d+b) - \sqrt{2k - (d-b)^2}\right) 
\]
and
\[
u(b,d,k) = \frac{1}{2}\left(-(d+b) + \sqrt{2k - (d-b)^2}\right).
\]
By analogy with the proof of Theorem~\ref{thm:nonlinearchance}, it suffices to show that $l$ is convex and $u$ is concave. To prove this, it suffices to show that $\sqrt{2k - (d-b)^2}$ is concave, which holds since $\sqrt{\cdot}$ is concave increasing and $2k-(d-b)^2$ is concave. Note we allow $\epsilon > \frac{1}{2}$ because this proof requires only monotonicity properties of $S_\epsilon$, not $\bar S_\epsilon$.
\end{proof}
\end{lemma}

We also observe that when $f$ and $g$ in Theorem~\ref{thm:nonlinearchance} are piecewise linear, e.g., as in Lemma~\ref{eq:absconvex}, then we have demonstrated the convexity of a special family of joint linear chance constraints.

\subsection{Distributionally robust two-sided chance constraints}
\label{sec:distrobust}

So far we have left unquestioned the assumption that the parameters $\mu$ and $\Sigma$ of the Gaussian distribution are known with certainty, when often they are subject to measurement error. For the case of linear chance constraints, Bienstock et al.~\cite{ccopf-sirev} propose a tractable model that enforces robustness with respect to deviations of the parameters $\mu$ and $\Sigma$ within a known uncertainty set $U$. Lubin et al.~\cite{JuMPChanceCaseStudy} implement this model and demonstrate significant cost savings in the context of short-term operational planning of power systems when tested against out-of-sample realizations of uncertainty. Here, we define and demonstrate tractability of a similar distributionally robust model in the context of two-sided chance constraints.

Let $\xi \sim N(\mu,\Sigma)$ be a jointly distributed Gaussian random vector with mean $\mu$ and positive definite covariance matrix $\Sigma$ and $0 < \epsilon \le \frac{1}{2}$, and let $LL^T = \Sigma$ be the Cholesky decomposition of $\Sigma$.

From Lemma~\ref{lem:generalcc}, recall
\begin{equation}
\mathbb{P}(a \le x^T\xi \le b) \geq 1 - \epsilon
\end{equation}
iff
\begin{equation}
\exists\, t \ge ||L^T x||_2 \text{ such that } (a-\mu^T x,b - \mu^T x,t) \in \bar S_\epsilon.
\end{equation}

We define the \textit{distributionally robust} (or \textit{ambiguous}) two-sided chance constraint as:
\begin{equation}\label{eq:robustcc}
\mathbb{P}_{\xi \sim N(\mu,\Sigma)}(a \le x^T\xi \leq b) \geq 1 - \epsilon \quad \forall (\mu,\Sigma) \in U
\end{equation}

\begin{lemma}
For $\epsilon \leq \frac{1}{2}$ and under the assumption that the uncertainty set decomposes by $\mu$ and $\Sigma$, i.e., $U = U_\mu \times U_\Sigma$, then the constraint~\eqref{eq:robustcc} is tractable if we can tractably optimize a linear objective over the sets $U_\mu$ and $U_\Sigma$.
\begin{proof}
Note that~\eqref{eq:robustcc} is a convex constraint, because it is the intersection of (infinitely) many convex constraints. We will prove tractability by demonstrating that we can easily separate, i.e., find the worst-case $\mu$ and $\Sigma$ given $(a,b,c)$.

We have that \eqref{eq:robustcc} holds iff $\exists\, t$ s.t.
\begin{align}
t &\geq ||L^T_{\Sigma}x|| &\forall \Sigma \in U_\Sigma \label{eq:conicrobust}\\
(a-\mu^Tx,b-\mu^Tx,t) &\in \bar S_\epsilon &\forall \mu \in U_\mu \label{eq:cdfrobust}
\end{align}
Constraint~\eqref{eq:conicrobust} can be reformulated as $t \geq \sqrt{\max_{\Sigma \in U_\Sigma}x^T\Sigma x}$, so we can separate by optimizing a linear objective over $U_\Sigma$. For $t>0$, the separation problem corresponding to constraint~\eqref{eq:cdfrobust} is
\begin{equation}\label{eq:cdfseparate}
\min_{\mu\in U_\mu}\Phi((b-\mu^Tx)/t) - \Phi((a-\mu^Tx)/t),
\end{equation}
which is a minimization of a log-concave function. However, observe that it is essentially a one-dimensional problem depending on $\mu^Tx$, so it can be solved by testing with the values $\min_{\mu \in U_\mu} \mu^Tx$ and $\max_{\mu \in U_\mu} \mu^Tx$.
\end{proof}
\end{lemma}

\section{Computational tractability of $S_\epsilon$ and $\bar S_\epsilon$}\label{sec:tractability}

We have demonstrated applications of the set $S_\epsilon$ and its conic hull $\bar S_\epsilon$ to represent a number of classes of convex chance constraints. So far, we have used these sets as a theoretical tool in order to prove convexity. In practice, we are also interested in computationally tractable representations of these sets, which ideally can be used within off-the-shelf solvers. 

\subsection{Representation using convex functions}\label{sec:nlp_representation}

Recall from Lemma~\ref{lem:seps} that the function $f(x,y) = \int_x^y \varphi(t)\,dt = \Phi(y)-\Phi(x)$ is log-concave. Therefore the equivalence $(x,y) \in S_\epsilon$ iff $\log f(x,y) \ge \log(1-\epsilon)$ provides a representation which can essentially be used directly within derivative-based nonlinear solvers, which expect constraints in the form $f(x) \le 0$ where $f$ is smooth and convex. Furthermore, the \textit{perspective function}
\[
    g(x,y,z) = z(\log(\Phi(y/z)-\Phi(x/z))-\log(1-\epsilon))
\]
is concave~\cite{HiriartLemarechal93book2}, and one can see that $(x,y,z) \in \bar S_\epsilon$ iff $g(x,y,z) \ge 0$ and $z \ge 0$, which provides a potentially useful representation of $\bar S_\epsilon$. 

The above representations are valid for the interval $\epsilon \in (0,1)$. For the special case of $\epsilon \in (0,\frac{1}{2})$, we note that for $x > 0$, $\Phi(x)$ is concave, and for $x < 0$, $\Phi(x)$ is convex. Since $\epsilon < \frac{1}{2}$ and $f(x,y) = \Phi(y)-\Phi(x) \ge 1-\epsilon$ imply $x < 0$ and $y > 0$, we note that $f$ itself is concave over the domain of $S_\epsilon$ because it is a sum of two concave functions. This observation provides an alternative convex representation of $(x,y,z) \in \bar S_\epsilon$ with the constraints
\begin{equation}\label{eq:convexfunc}
    z(\Phi(y/z)-\Phi(x/z) - (1-\epsilon)) \ge 0 \text{ and } z \ge 0.
\end{equation}
With either convex representation, the derivatives are easy to compute when $z > 0$. However, derivative-based solvers may fail as $z \to 0$.

\subsection{Separation oracles}

A functional, derivative-based representation of $\bar S_\epsilon$ may be directly applicable in many situations, but alternative solution methods exist. For example, algorithms for convex mixed-integer nonlinear optimization typically make use of a combination of continuous nonlinear relaxations and iteratively generated polyhedral outer approximations~\cite{StubbsMehrotra99,Bonmin}. In this section we discuss how separation oracles could be implemented to generate such polyhedral outer approximations. Our focus is on developing separation oracles which lead to polyhedral approximations which are ``better'' than the more common approach which follows from the functional representation of Section~\ref{sec:nlp_representation}, in the sense of producing hyperplanes which are tangent to the set $\bar S_\epsilon$.

In brief, when a separation oracle for the set $S_\epsilon$ is given a point $(x,y)$, it first determines if $(x,y) \in S_\epsilon$. If $(x,y) \not \in S_\epsilon$, it returns a hyperplane $(a,b) \in \mathbb{R}^2 \times \mathbb{R}$ such that $a_1x + a_2y > b$ and $S_\epsilon$ is contained in the halfspace defined by $\{ (x,y) : a_1 x + a_2 y \le b\}$. Hence, the hyperplane separates the point $(x,y)$ from the set $S_\epsilon$. 

First we note that a separation oracle for $S_\epsilon$ immediately provides a separation oracle for the conic hull $\bar S_\epsilon$. Suppose $(x,y,z) \not \in \bar S_\epsilon$ and $z > 0$. Then $(x/z,y/z) \not \in S_\epsilon$ so we take a hyperplane $a_1 x + a_2 y = b$ which separates $(x/z,y/z)$ from $S_\epsilon$, then the hyperplane $a_1 x + a_2 y - b z = 0$ separates $(x,y,z)$ from $\bar S_\epsilon$. If $z = 0$, note that, assuming $\epsilon \le \frac{1}{2}$, $(x,y,0) \in \bar S_\epsilon$ iff $y \ge 0$ and $x \le 0$, so these two constraints serve as the separating hyperplanes in this case. Thus we restrict our discussion to separation oracles for $S_\epsilon$.

The most straightforward separation oracle for a convex set described by a smooth convex function is as follows.
For any smooth, convex function $f$, if we are given $x'$ with $f(x') > 0$, then the hyperplane $f(x') + \nabla f(x')(x-x') \le 0$ separates $x'$ from the feasible set of $\{ x : f(x) \le 0 \}$~\cite{Bonmin}\footnote{Taking for granted that we can compute $\Phi(\cdot)$ efficiently, this gradient-based separation oracle may also be of theoretical use in proving tractability via the ellipsoid algorithm~\cite{schrijver2003}.}. For the case of $S_\epsilon$, however, this hyperplane is weak. More specifically, we have
\begin{equation}\label{eq:sepnlp}
    (x,y) \in S_\epsilon \text{ iff } f(x,y) := 1-\epsilon-\Phi(y)+\Phi(x) \le 0.
\end{equation}

\begin{figure}[ht]
    \centering
    \includegraphics[scale=0.3]{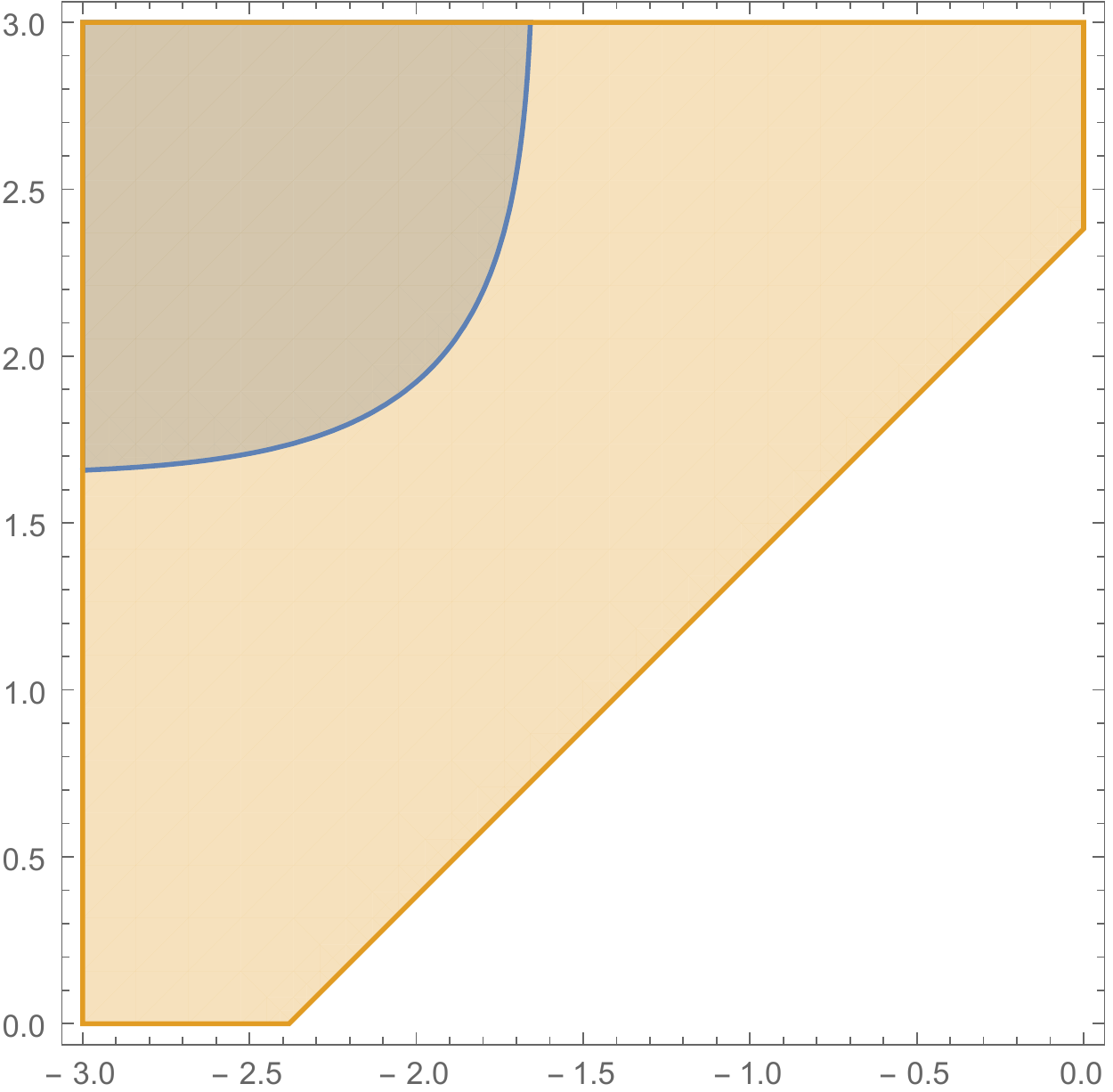}
    \caption{In blue, the set $S_\epsilon$, with $\epsilon=0.05$. In orange, the half space corresponding to the separating hyperplane~\eqref{eq:badhyperplane}. The hyperplane separates the point $(0,0)$ but is not tangent to $S_\epsilon$.}\label{fig:sepnlp}
\end{figure}

If we use this representation to separate the point $(0,0)$, then $f(0,0) = 1-\epsilon$ and $\nabla f(0,0) = (\frac{1}{\sqrt{2\pi}},-\frac{1}{\sqrt{2\pi}})$, and our separating hyperplane is 
\begin{equation}\label{eq:badhyperplane}
x-y \le -\sqrt{2\pi}(1-\epsilon).
\end{equation}
Figure~\ref{fig:sepnlp} shows the set $S_\epsilon$ together with this separating hyperplane. Observe that the hyperplane is not tangent to $S_\epsilon$, which means that it may serve poorly as an outer approximation.

Instead of using this hyperplane, we might consider computing the best possible separating hyperplane with the same slope by evaluating $\max_{(x,y) \in S_\epsilon} x - y$. By a symmetry argument which we omit here, this value is $2 \Phi^{-1}(\epsilon/2)$, so we can strengthen the previous hyperplane to
\begin{equation}\label{eq:betterhyperplane}
    x-y \le 2 \Phi^{-1}(\epsilon/2).
\end{equation}

More generally, the \textit{support function} $\sigma_{S_\epsilon}(a,b) = \max_{(x,y)\in S_\epsilon} ax+by$ enables one to compute the best possible separating hyperplane of a given slope~\cite{HiriartLemarechal93book2}. Another approach to generating tangent separating hyperplanes is to compute an orthogonal projection of the point $(x,y)$ onto the set $S_\epsilon$ and then add a hyperplane which is tangent to the projected point. In the following section, we provide simple representations of the support function and orthogonal projection operators for the set $S_\epsilon$. These developments may enable practical implementations.

\subsection{Support function of $S_\epsilon$ and orthogonal projection onto $S_\epsilon$}

We begin with a lemma which characterizes the boundary of $S_\epsilon$.

\begin{lemma}
    Let $\epsilon \in (0,\frac{1}{2})$. Then the point $(x,y)$ lies on the boundary of the set $S_\epsilon$, i.e., $\Phi(y)-\Phi(x) = 1-\epsilon$ iff $\exists\, \lambda \in (0,1)$ such that $x = \Phi^{-1}(\lambda\epsilon)$ and $y = \Phi^{-1}(1-(1-\lambda)\epsilon)$.
   \begin{proof}
       First, let $\lambda \in (0,1)$. Then $\Phi(\Phi^{-1}(1-(1-\lambda)\epsilon)) - \Phi(\Phi^{-1}(\lambda\epsilon)) = (1-(1-\lambda)\epsilon) - \lambda\epsilon = 1-\epsilon$, and so the point $(\Phi^{-1}(\lambda\epsilon),\Phi^{-1}(1-(1-\lambda)\epsilon))$ lies on the boundary. In the other direction, suppose the point $(x,y)$ is on the boundary of $S_\epsilon$. Set $\lambda = \Phi(x)/\epsilon$. Note that since $x$ and $y$ are finite, we must have
       $0 < \Phi(x) < \epsilon$ and $1-\epsilon < \Phi(y) < 1$, and hence $0 < \lambda < 1$. Trivially $x = \Phi^{-1}(\lambda\epsilon)$. Then $\Phi(y) = 1-\epsilon+\Phi(x)$, and we see that $y = \Phi^{-1}(1-\epsilon+\lambda\epsilon) = \Phi^{-1}(1-(1-\lambda)\epsilon)$. 
   \end{proof}
\end{lemma}

This result provides an explicit univariate parameterization of the boundary of $S_\epsilon$ in terms of $\lambda$, which is quite useful for computational purposes. For example, suppose we wanted to minimize a function $g(x,y)$ along the boundary of $S_\epsilon$. Then this problem can be formulated as a one-dimensional search problem,
\begin{equation}\label{eq:optboundary}
\min_{\lambda \in (0,1)} g(\Phi^{-1}(\lambda\epsilon), \Phi^{-1}(1-(1-\lambda)\epsilon)).
\end{equation}

The following lemma uses the formulation~\eqref{eq:optboundary} to demonstrate that optimization of some linear functions over $S_\epsilon$ can be expressed as a univariate convex optimization problem.

\begin{lemma}\label{lem:linear}
Suppose $g(x,y) = ax+by$ with $a < 0$ and $b > 0$. Then \eqref{eq:optboundary} is a smooth, strictly convex optimization problem.
\begin{proof}
Define $h(\lambda) := g(\Phi^{-1}(\lambda\epsilon), \Phi^{-1}(1-(1-\lambda)\epsilon))$. We explicitly calculate the derivatives using the following basic formulas:

\[
\frac{d}{dx} \Phi^{-1}(x) = \sqrt{2\pi} e^{\frac{\Phi^{-1}(x)^2}{2}}
\]
\[
\frac{d^2}{dx^2} \Phi^{-1}(x) = 2\pi\Phi^{-1}(x) e^{\Phi^{-1}(x)^2}
\]

so
\[
\frac{d^2}{d\lambda^2} h(\lambda) = 2\pi a\epsilon^2 \Phi^{-1}(\lambda\epsilon) e^{\Phi^{-1}(\lambda\epsilon)^2} + 2\pi b\epsilon^2\Phi^{-1}(1-(1-\lambda)\epsilon) e^{\Phi^{-1}(1-(1-\lambda)\epsilon)^2}.
\]

Note $\Phi^{-1}(\lambda\epsilon) < 0$ and $\Phi^{-1}(1-(1-\lambda)\epsilon) > 0$, so given $a < 0$ and $b > 0$, we have that $\frac{d^2}{d\lambda^2} h(\lambda) > 0$.
\end{proof}
\end{lemma}

Following Lemma~\ref{lem:linear} we have an efficient way to evaluate the support function
\[
\sigma_{S_\epsilon}(a,b) = \max_{(x,y)\in S_\epsilon} ax+by.
\]
Specifically, when $a > 0$ and $b < 0$, we solve a one-dimensional convex minimization problem. If $a < 0$ or $b > 0$, then $\sigma_{S_\epsilon}(a,b) = \infty$. If $a = 0$ and $b < 0$, $\sigma_{S_\epsilon}(a,b) = b\Phi(1-\epsilon)$. If $b = 0$ and $a > 0$, $\sigma_{S_\epsilon}(a,b) = a\Phi(\epsilon)$. These last two cases follow from taking the limit when $\lambda = 0$ and $\lambda = 1$, respectively.

\begin{lemma}
We can compute an orthogonal projection onto $S_\epsilon$ by solving a one-dimensional strictly convex minimization problem.
\begin{proof}
    Similar to Lemma~\ref{lem:linear}, we will use the parameterization of the boundary, but solving~\eqref{eq:optboundary} over a restricted domain. Given $(a,b) \not \in S_\epsilon$, the orthogonal projection is the solution to~\eqref{eq:optboundary} with $g(x,y) = \frac{1}{2}(x-a)^2+\frac{1}{2}(y-b)^2$. Actually we do not need to optimize over all $\lambda \in (0,1)$; note that the orthogonal projection always lies on the boundary of $S_\epsilon$ between the projections along the x and y axes. More specifically, we need only consider
\begin{equation}\label{eq:proj_domain}
\lambda \in \left(1-\frac{1}{\epsilon}(1-\Phi(b)), \frac{1}{\epsilon}\Phi(a)\right),
\end{equation}
and within this interval, by construction, the inequalities 
\begin{equation}\label{eq:proj_ineq}
\Phi^{-1}(\lambda\epsilon) \leq a \text{ and }\Phi^{-1}(1-(1-\lambda)\epsilon) \geq b
\end{equation}
hold.

Define $h(\lambda) := \frac{1}{2}(\Phi^{-1}(\lambda\epsilon)-a)^2+\frac{1}{2}(\Phi^{-1}(1-(1-\lambda)\epsilon)-b)^2$. We will prove strict convexity of $h$ within the domain~\eqref{eq:proj_domain} by showing that $\frac{d^2h}{d\lambda^2} > 0$.  From the chain rule (for arbitrary $f$),
\[
\frac{d^2}{dx^2} \frac{1}{2}(f(x)-a)^2 = \left(\frac{df}{dx}(x)\right)^2 + (f(x)-a)\frac{d^2f}{dx^2}(x).
\]
Discarding the squared first derivative terms, we have 
\[
\frac{d^2h}{d\lambda^2}(\lambda) \geq (\Phi^{-1}(\lambda\epsilon)-a) \frac{d^2}{d\lambda^2} \Phi^{-1}(\lambda\epsilon) +
(\Phi^{-1}(1-(1-\lambda)\epsilon)-b)  \frac{d^2}{d\lambda^2}(\Phi^{-1}(1-(1-\lambda)\epsilon)
\]
The result follows from noting that $\frac{d^2}{d\lambda^2} \Phi^{-1}(\lambda\epsilon) < 0$ and $\frac{d^2}{d\lambda^2}(\Phi^{-1}(1-(1-\lambda)\epsilon) > 0$ combined with the inequalities~\eqref{eq:proj_ineq}.
\end{proof}
\end{lemma}

\subsection{An approximate polyhedral representation of $S_\epsilon$}\label{sec:polyapprox}

In this section, we develop an approximate polyhedral representation of $S_\epsilon$.

\begin{definition}
    A polyhedron $P_\epsilon$ is an \textit{outer approximation} of $S_\epsilon$ if $S_\epsilon \subset P_\epsilon$.
\end{definition}

While polyhedral outer approximations are straightforward to generate, either through an iterative cutting-plane procedure or by preselecting a number of tangent hyperplanes, we are interested in outer approximations with a provable approximation guarantee, in the sense which we now define.

\begin{definition}
    A family of polyhedral outer approximations $P_\epsilon$ forms an $\alpha$-approximation of $S_\epsilon$ if $\forall \epsilon \in (0,\frac{1}{2}]$,
\begin{equation}
    \Phi(y)-\Phi(x) \ge 1-\alpha\epsilon\quad \forall (x,y) \in P_\epsilon.
\end{equation}
Or equivalently, when $\alpha\epsilon < 1$, $S_\epsilon \subset P_\epsilon \subset S_{\alpha\epsilon}$.
\end{definition}

We restrict $\epsilon \le \frac{1}{2}$ for notational convenience and because this is the case of direct interest, although many of the results here generalize for $\epsilon \in (0,1)$.

Note that although our development is from the perspective of outer approximation, a family of polyhedral outer approximations may be used to generate conservative approximations as well, since if $P_\epsilon$ is an $\alpha$-approximation, then
\[
    P_{\epsilon/\alpha} \subset S_\epsilon\quad \forall \epsilon \in (0,\frac{1}{2}].
\]

We begin with a very simple 2-approximation of $S_\epsilon$ with the axis-aligned polyhedra:
\begin{equation}
    A_\epsilon = \{ (x,y) \in \mathbb{R}^2 : x \le \Phi^{-1}(\epsilon), y \ge \Phi^{-1}(1-\epsilon) \}.
\end{equation}
For $(x,y) \in A_\epsilon$, by monotonicity of the cumulative density function $\Phi$ we conclude
\begin{equation}
    \Phi(y) - \Phi(x) \ge (1-\epsilon) - \epsilon \ge 1-2\epsilon.
\end{equation}

This 2-approximation of $S_\epsilon$ is equivalent to representing $\mathbb{P}(a \le x^T\xi \le b) \ge 1-\epsilon$ by using the two standard linear chance constraints $\mathbb{P}(a \le x^T\xi) \ge 1-\epsilon$ and $\mathbb{P}(x^T\xi \le b) \ge 1-\epsilon$. Bienstock et al.~\cite{ccopf-sirev} employ this approximation citing improved computational tractability.

The 2-approximation model is the best one can achieve with two linear constraints in the following sense. The set $S_\epsilon$ has two extreme rays: $(-1,0)$ and $(0,1)$, which follow from the fact that $\Phi$ is monotonic increasing. Therefore, any outer approximation of $S_\epsilon$ must contain these rays. If, in addition, these are not the extreme rays of the outer approximation, then the approximation \textit{cannot} be an $\alpha$-approximation for any $\alpha$, because $S_{\alpha\epsilon}$ cannot contain the set.

\begin{figure}
\centering
\includegraphics[scale=0.4]{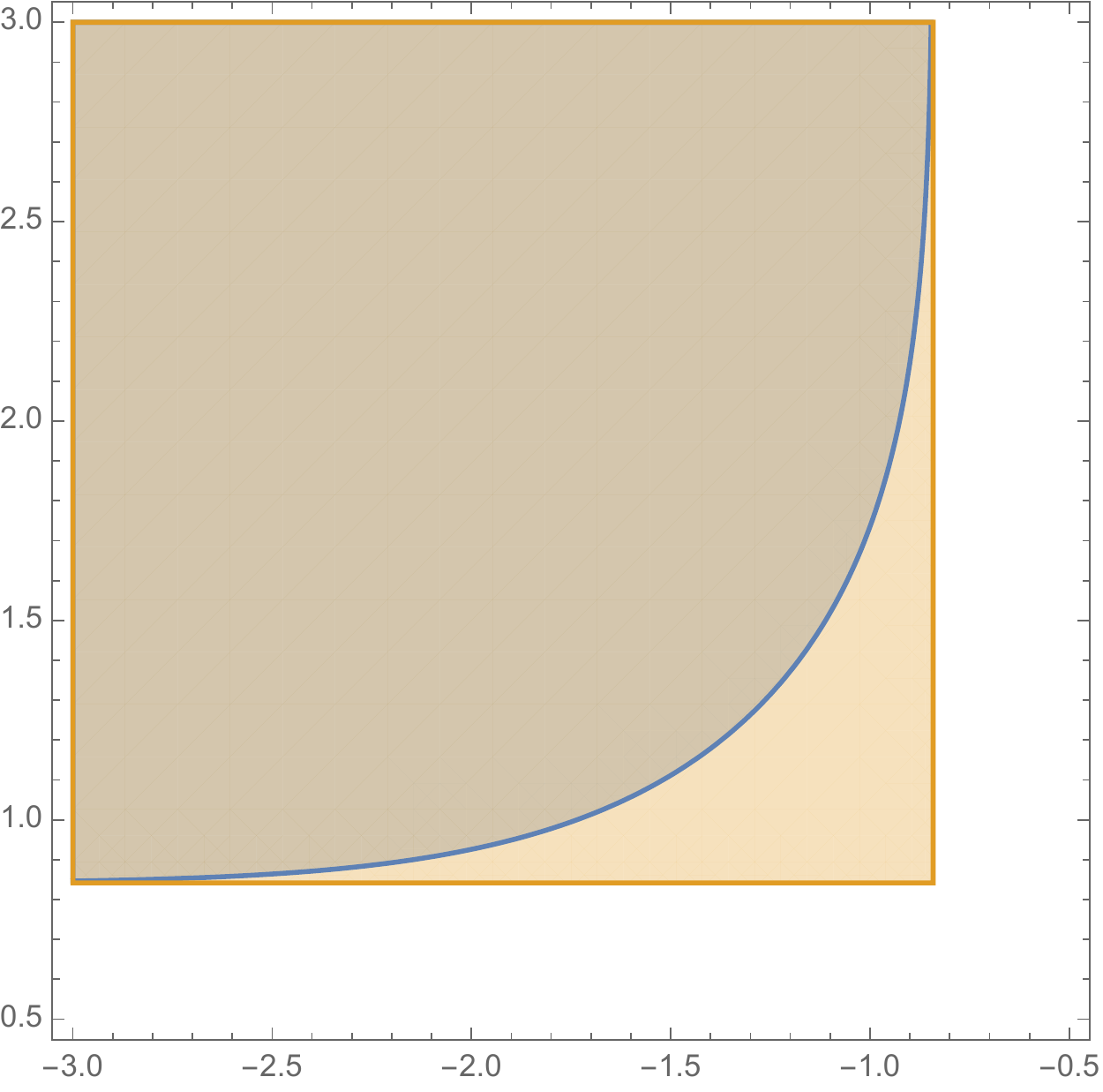}
\hspace{1cm}
\includegraphics[scale=0.4]{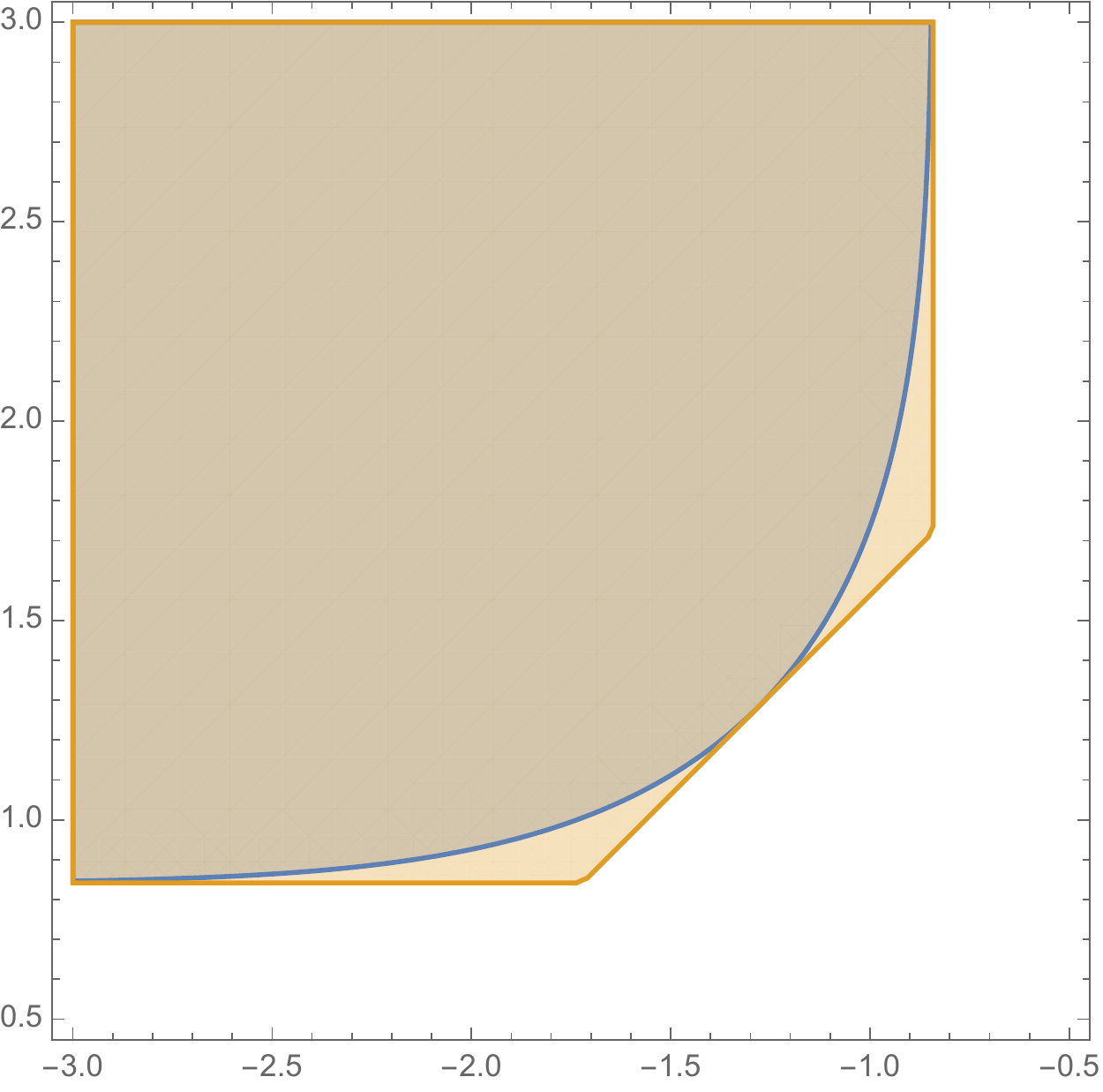}
\caption{In blue, the set $S_\epsilon$. In orange, the polyhedral outer approximation $A_\epsilon$ (left) and $B_\epsilon$ (right). By adding a single additional inequality, we strengthen the relaxation significantly.}\label{fig:twothreecut}
    \end{figure}

The main result of this section is that with a single additional linear constraint, one may improve the above 2-approximation to a 1.25-approximation. The axis-aligned approximation performs poorly at the ``corner'' where $x = \Phi^{-1}(\epsilon)$ and $y = \Phi^{-1}(1-\epsilon)$. If we add a hyperplane to separate this point, from the previous discussion we obtain the hyperplane~\eqref{eq:betterhyperplane}.

Therefore we define the family of polyhedra as
\begin{equation}\label{eq:3cutdefn}
    B_\epsilon := \{ (x,y) \in \mathbb{R}^2 : x \le \Phi^{-1}(\epsilon), y \ge \Phi^{-1}(1-\epsilon),  x-y \le 2 \Phi^{-1}(\epsilon/2) \}.
\end{equation}
The family $B_\epsilon$ forms a valid outer approximation because $A_\epsilon$ is a valid family, and we've added a valid separating hyperplane. Figure~\ref{fig:twothreecut} displays the two families of approximations for a fixed $\epsilon$.

The following lemma simplifies the task of proving the $\alpha$-approximation.

\begin{lemma}
    For $\alpha < 2$, a family of polyhedral outer approximations $P_\epsilon$ forms an $\alpha$-approximation of $S_\epsilon$ iff $\forall \epsilon \in (0,\frac{1}{2}]$
\begin{enumerate}
    \item $\forall$ vertices $(x,y)$ of $P_\epsilon,$ we have $\Phi(y) - \Phi(x) \ge 1-\alpha\epsilon$, and
    \item the extreme rays of $P_\epsilon$ are $(-1,0)$ and $(0,1)$.
\end{enumerate}
That is, it is sufficient to verify the approximation quality at the vertices.
\begin{proof}
    Fix $\epsilon$ and suppose that the two above conditions hold. Then all vertices, by definition, are contained in the set $S_{\alpha\epsilon}$ (our assumptions imply $\alpha\epsilon < 1$). By convexity of $S_{\alpha\epsilon}$, this implies that all convex combinations of the vertices of $P_\epsilon$ are contained in $S_{\alpha\epsilon}$. All elements of the polyhedron $P_\epsilon$ can be represented as a convex combination of its vertices plus a conic combination of its extreme rays. Since the extreme rays $(-1,0)$ and $(0,1)$ are also extreme rays of $S_{\alpha\epsilon}$, it follows that $P_{\epsilon} \subset S_{\alpha\epsilon}$. 
\end{proof}
\end{lemma}

The vertices of $B_\epsilon$ are $(\Phi^{-1}(\epsilon),\Phi^{-1}(\epsilon)-2\Phi^{-1}(\epsilon/2))$ and $(2\Phi^{-1}(\epsilon/2)+\Phi^{-1}(1-\epsilon),\Phi^{-1}(1-\epsilon))$. By the identities $\Phi^{-1}(1-\epsilon) = -\Phi^{-1}(\epsilon)$ and $\Phi(-x) = 1-\Phi(x)$ we see that these vertices are in fact symmetric, so it is sufficient to consider only one of them.

We first establish a simple bound that does not use any deep properties of the Gaussian distribution.

\begin{lemma}
The ``three-cut'' family of outer approximations $B_\epsilon$ forms a 1.5-approximation of $S_\epsilon$.
\begin{proof}
Consider the vertex $(\Phi^{-1}(\epsilon),\Phi^{-1}(\epsilon)-2\Phi^{-1}(\epsilon/2))$. It is sufficient to show that it is contained in the set $S_{1.5\epsilon}$.
\begin{align}
    \Phi(\Phi^{-1}(\epsilon)-2\Phi^{-1}(\epsilon/2)) - \Phi(\Phi^{-1}(\epsilon)) &=
    \Phi((\Phi^{-1}(\epsilon)-\Phi^{-1}(\epsilon/2))-\Phi^{-1}(\epsilon/2)) - \epsilon\\
    &> \Phi(-\Phi^{-1}(\epsilon/2)) - \epsilon \\
    &= \Phi(\Phi^{-1}(1-\epsilon/2)) - \epsilon\\
    &= 1-1.5\epsilon,
\end{align}
where the inequality follows from $\Phi^{-1}(\epsilon) > \Phi^{-1}(\epsilon/2)$ and monotonicity of $\Phi$.
\end{proof}
\end{lemma}

We can improve this bound by using properties of the Gaussian distribution. Lemmas~\ref{lem:tail2} and \ref{lem:4bound} below develop the necessary properties, and Theorem~\ref{thm:125} states the final result.

\begin{lemma}\label{lem:tail2}
$\Phi^{-1}(1-\frac{\epsilon}{2}) - \Phi^{-1}(1-\epsilon) \geq \Phi^{-1}(1-\frac{\epsilon}{4}) - \Phi^{-1}(1-\frac{\epsilon}{2})$ for $\epsilon \in (0,\frac{1}{2}]$.
\begin{proof}
Let $f(\epsilon) = \Phi^{-1}(1-\frac{\epsilon}{2}) - \Phi^{-1}(1-\epsilon)$. Then we intend to show $f(\epsilon) \geq f(\epsilon/2)\, \forall \epsilon \in (0,\frac{1}{2}]$. It suffices to show that $f$ is monotonic increasing over the interval.

Recalling
\[
\frac{d}{dx} \Phi^{-1}(x) = \sqrt{2\pi} \exp\left(\frac{\Phi^{-1}(x)^2}{2}\right),
\]
we have
\[
f'(\epsilon)= -\frac{1}{2}\sqrt{2\pi}\exp\left(\frac{\Phi^{-1}(1-\frac{\epsilon}{2})^2}{2}\right)+\sqrt{2\pi}\exp\left(\frac{\Phi^{-1}(1-\epsilon)^2}{2}\right).
\]
We will show that $f'$ is always positive for $\epsilon \in (0,\frac{1}{2}]$. At $\epsilon = \frac{1}{2}$, $\Phi^{-1}(\frac{1}{2}) = 0$, so
\[
f'\left(\frac{1}{2}\right) = \sqrt{2\pi} - \frac{1}{2}\sqrt{2\pi} \exp\left(\Phi^{-1}\left(\frac{3}{4}\right)^2/2\right) \approx 0.93 > 0.
\]
Suppose, for contradiction, $f'(\epsilon') = 0$ for some $\epsilon'$. Then
\[
\exp\left(\frac{\Phi^{-1}(1-\epsilon')^2}{2}\right) = \frac{1}{2}\exp\left(\frac{\Phi^{-1}(1-\frac{\epsilon'}{2})^2}{2}\right)
\]
which implies
\begin{equation}\label{eq:fprimezero}
\Phi^{-1}(1-\epsilon')^2 = -2\log(2) + \Phi^{-1}\left(1-\frac{\epsilon'}{2}\right)^2.
\end{equation}

Note that $g(\epsilon) := \Phi^{-1}(1-\epsilon)^2$ is strictly convex for $\epsilon \in (0,1)$  by examination of the second derivative. This means that $g'(\epsilon)$ is strictly monotonic increasing. We're looking for a solution to $g(\epsilon/2) - g(\epsilon) = 2\log(2)$. Note that $g(\epsilon/2) - g(\epsilon)$ is strictly decreasing over the interval because $(1/2)g'(\epsilon/2) - g'(\epsilon) < 0$.  One can verify the limit
\[
\lim_{\epsilon\to 0+} \Phi^{-1}\left(1-\frac{\epsilon}{2}\right)^2 -\Phi^{-1}(1-\epsilon)^2  = 2\log(2),
\]
which implies in fact that there can be no solution to~\eqref{eq:fprimezero}. This proves our original claim.

\end{proof}
\end{lemma}

\begin{lemma}\label{lem:4bound}
$\Phi^{-1}(\epsilon) - 2\Phi^{-1}(\frac{\epsilon}{2}) \geq \Phi^{-1}(1-\frac{\epsilon}{4})$ for $\epsilon \in (0,\frac{1}{2}]$
\begin{proof}
Applying Lemma~\ref{lem:tail2}, we have:
\begin{align}
\Phi^{-1}(\epsilon) - 2\Phi^{-1}(\frac{\epsilon}{2}) &= \Phi^{-1}(1-\frac{\epsilon}{2}) + (\Phi^{-1}(1-\frac{\epsilon}{2}) - \Phi^{-1}(1-\epsilon))\\
&\geq \Phi^{-1}(1-\frac{\epsilon}{2}) + (\Phi^{-1}(1-\frac{\epsilon}{4}) - \Phi^{-1}(1-\frac{\epsilon}{2}))\\
&=\Phi^{-1}(1-\frac{\epsilon}{4})
\end{align}
\end{proof}
\end{lemma}

\begin{theorem}\label{thm:125}
The ``three-cut'' family of outer approximations $B_\epsilon$ forms a 1.25-approximation of $S_\epsilon$.
\begin{proof}
Consider the vertex $(\Phi^{-1}(\epsilon),\Phi^{-1}(\epsilon)-2\Phi^{-1}(\epsilon/2))$. It is sufficient to show that it is contained in the set $S_{1.25\epsilon}$.
\begin{align}
    \Phi(\Phi^{-1}(\epsilon)-2\Phi^{-1}(\epsilon/2)) - \Phi(\Phi^{-1}(\epsilon)) &=\\
    &\ge \Phi(\Phi^{-1}(1-\epsilon/4)) - \epsilon \\
    &= 1-1.25\epsilon,
\end{align}
where the inequality follows from Lemma~\ref{lem:4bound}.
\end{proof}
\end{theorem}

With an additionally highly technical argument which we omit for brevity, it is possible to show that the 1.25 value is tight; that is, the ``three-cut`` family of outer approximations $B_\epsilon$ is \textit{not} an $\alpha$-approximation for any $\alpha < 1.25$.

We summarize the results of this section with a succinct statement of an SOC outer approximation of the two-sided chance constraint based on $B_\epsilon$.

\begin{lemma}\label{lem:generalcc_1.25}
Let $\xi \sim N(\mu,\Sigma)$ be a jointly distributed Gaussian random vector with mean $\mu$ and positive definite covariance matrix $\Sigma$ and $0 < \epsilon \le \frac{1}{2}$. Let $LL^T = \Sigma$ be the Cholesky decomposition of $\Sigma$.
The following extended formulation, with the additional variable $t$,
\begin{align}
    t \ge & ||L^Tx||_2,\label{eq:outer1}\\
    a-\mu^Tx & \le \Phi^{-1}(\epsilon)t,\label{eq:outer2}\\
    b - \mu^T x & \ge \Phi^{-1}(1-\epsilon)t,\label{eq:outer3}\\
    a-b & \le  2 \Phi^{-1}(\epsilon/2)t.\label{eq:outer4}
\end{align}
is an SOC outer approximation of the constraint
\[
\mathbb{P}(a \le x^T\xi \leq b) \geq 1 - \epsilon
\]
which in fact guarantees
\[
\mathbb{P}(a \le x^T\xi \leq b) \geq 1 - 1.25\epsilon.
\]
\begin{proof}
From Lemma~\ref{lem:generalcc},
\begin{equation}
\mathbb{P}(a \leq x^T\xi \leq b) \geq 1-\epsilon
\end{equation}
iff
\begin{equation}
\exists\, t \ge ||L^T x||_2 \text{ such that } (a-\mu^T x,b - \mu^T x,t) \in \bar S_\epsilon.
\end{equation}
We take the conic hull of the polyhedral representation $B_\epsilon$~\eqref{eq:3cutdefn} of $S_\epsilon$ in order to represent $\bar S_\epsilon$.
\end{proof}
\end{lemma}

\section{Approximation of quadratic chance constraints}\label{sec:quadapprox}

Having extensively discussed the tractability of the two-sided chance constraint model and its extensions to represent more complex nonlinear chance constraints exactly, we return to our original motivation as discussed in Section~\ref{sec:motivation}. In this section, we will investigate the use of two-sided chance constraints to \textit{approximately} represent a family of challenging \textit{quadratic chance constraints}. These sets are of the form,
\begin{equation}\label{eq:quadprob}
H_\epsilon = \left\{ (a,b,c,d,k) \in \mathbb{R}^n\times \mathbb{R} \times \mathbb{R}^n \times \mathbb{R} \times \mathbb{R} : \mathbb{P}((a^T\xi + b)^2 + (c^T\xi + d)^2 \le k) \ge 1-\epsilon \right\},
\end{equation}
where $a$ and $c$, and $b$, $d$, and $k$ are (vector and scalar, resp.) decision variables and $\xi$ follows a multivariate Gaussian distribution with known mean and covariance matrix.

\subsection{Convexity of the quadratic chance constraint}

We are unaware of any existing results on the convexity of the set $H_\epsilon$~\eqref{eq:quadprob}. We present here a proof of nonconvexity for the case of $\epsilon = 0.455$. The counterexample, while not as strong as a proof of nonconvexity for \textit{all} $\epsilon \in (0,\frac{1}{2}]$, suggests that convexity is, at the least, not a simple extension of existing results such as those presented in Section~\ref{sec:exact_extensions} which hold for all $\epsilon \in (0,\frac{1}{2}]$. We leave the question of convexity of $H_\epsilon$ over the full range of $\epsilon$ for future work. Nevertheless, we take this counterexample as a justification for seeking tractable, convex approximations of $H_\epsilon$ in subsequent sections.

Consider the constraint
\begin{equation}\label{eq:nonconvex_example}
\mathbb{P}((x\xi_1)^2 + (y\xi_2)^2 \le 1) \ge 1-\epsilon,
\end{equation}
where $\xi_1$ and $\xi_2$ are independent, standard Gaussian random variables.
The constraint~\eqref{eq:nonconvex_example} is a special case of~\eqref{eq:quadprob} with $\xi = (\xi_1,\xi_2)$, $a = (x,0), b = 0, c = (0,y), d = 0$, and $k = 1$.

\begin{figure}[ht]
\centering
\includegraphics[scale=0.7]{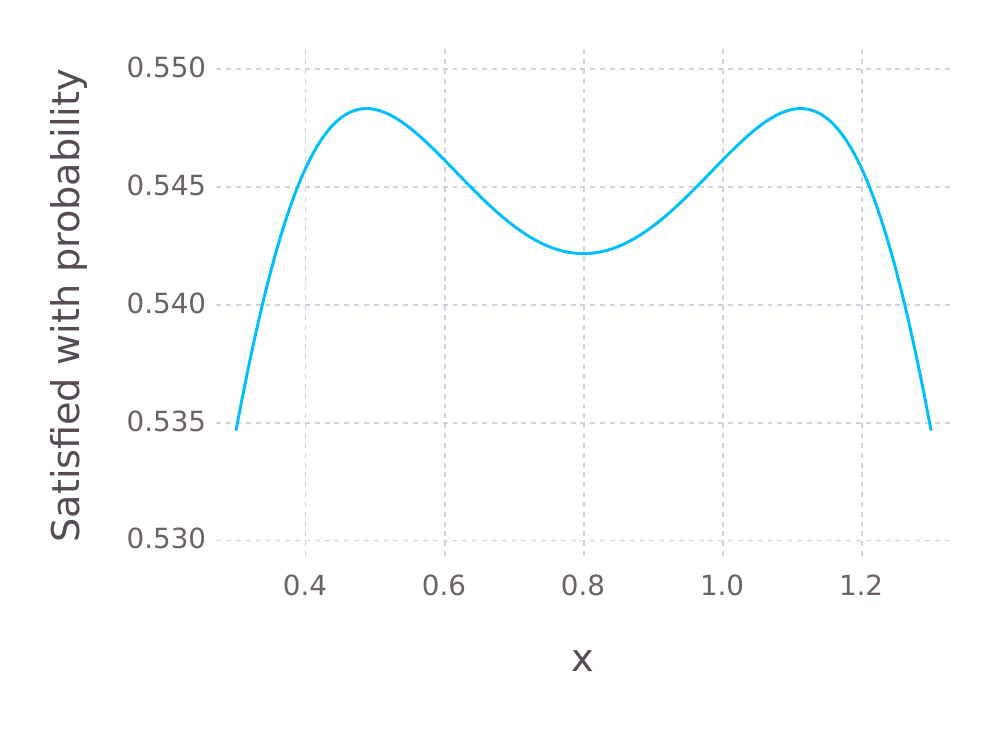}
\caption{On the vertical axis, the value of the left-hand side of~\eqref{eq:nonconvex_example} evaluated at the point $(x,-x+1.6)$ by numerical integration (with approximate error bounds of $10^{-7}$). We see that the set of points along this line that satisfy the quadratic chance constraint with probability 0.545 or greater, for example, is not convex. This proves nonconvexity of $H_\epsilon$ with $\epsilon = 0.455$. }
\label{fig:nonconvex_example}
\end{figure}

Figure~\ref{fig:nonconvex_example} traces the value of the left-hand side of~\eqref{eq:nonconvex_example} along the line $y = -x + 1.6$. We see that the upper level sets of the function $f(x,y) = \mathbb{P}((x\xi_1)^2 + (y\xi_2)^2 \le 1)$ are not convex. In particular, the points $(0.6,1.0)$ and $(1.0,0.6)$ belong to $H_{0.455}$ (by numerical integration with reported error bounds of $10^{-7}$) but the point $(0.8,0.8)$, their average, does not. We can evaluate $\mathbb{P}((0.8\xi_1)^2 + (0.8\xi_2)^2 \le 1)$ more explicitly as $F_{\chi_2}(1/0.8) \approx 0.542$ where $F_{\chi_2}$ is the cumulative distribution function of the chi distribution with two degrees of freedom.

\subsection{Approximation using two-sided constraints}

We propose an approximation of the quadratic chance constraint \eqref{eq:quadprob} by two absolute value constraints, essentially splitting up the squared terms into separate constraints. We use the union bound to enforce a conservative approximation, hence we introduce a parameter $\beta \in (0,1)$ to balance the trade-off between violations in the two separate constraints. In Lemma~\ref{lem:quadapprox} below, we state this formulation formally and prove that it is a valid convex, conservative approximation of \eqref{eq:quadprob} in an extended set of variables.

\begin{lemma}\label{lem:quadapprox}
Fix $\epsilon < \frac{1}{2}$, fix $\beta \in (0,1)$, and let
\begin{subequations}\label{eq:quadapproxall}
\begin{align}
G_{\epsilon,\beta} = \biggl\{ (a,b,c,d,k, f_1,f_2) \in \mathbb{R}^{2n+5} :& 
\mathbb{P}(|a^T\xi + b| \le f_1) \ge 1-\beta\epsilon \label{eq:abs1}\\
&\mathbb{P}(|c^T\xi + d| \le f_2) \ge 1-(1-\beta)\epsilon \label{eq:abs2} \\
&f_1^2 + f_2^2 \le k & \biggr\}.\label{eq:quadquad}
\end{align}
\end{subequations}
Let $G_{\epsilon,\beta}^{proj}$ be the projection of the set $G_{\epsilon,\beta}$ onto the variables $(a,b,c,d,k)$. Then $G_{\epsilon,\beta}^{proj}$ is convex and $G_{\epsilon,\beta}^{proj} \subseteq H_\epsilon$. That is, the set $G_{\epsilon,\beta}^{proj}$ is a conservative, convex approximation of the quadratic chance constraint~\eqref{eq:quadprob}.

\begin{proof}
Let $(a,b,c,d,k,f_1,f_2) \in G_{\epsilon,\beta}$. To simplify the proof, let $\chi = a^T\xi + b$ and $\psi = c^T\xi + d$ be random variables. Then $\chi^2 \le f_1^2$ and $\psi^2 \le f_2^2$ implies $\chi^2 + \psi^2 \le k$, which gives us the inequality
\begin{equation}
\mathbb{P}(\chi^2 + \psi^2 \le k) \ge \mathbb{P}(\chi^2 \le f_1^2 \text{ and } \psi^2 \le f_2^2).
\end{equation}
From the union bound,
\begin{align}
\mathbb{P}(\chi^2 \le f_1^2 \text{ and } \psi^2 \le f_2^2) &\ge \mathbb{P}(\chi^2 \le f_1^2) + \mathbb{P}(\psi^2 \le f_2^2) - 1 \\
&\ge 1-\beta\epsilon + 1 - (1-\beta)\epsilon - 1\\
&=1 -\epsilon,
\end{align}
which proves $G_{\epsilon,\beta}^{proj} \subseteq H_\epsilon$. Convexity of $G_{\epsilon,\beta}$ (and therefore $G_{\epsilon,\beta}^{proj}$) follows from Lemma~\ref{eq:absconvex} and the fact that $f_1^2 + f_2^2 \le k$ is a convex quadratic constraint.
\end{proof}
\end{lemma}

Replacing constraints~\eqref{eq:abs1} and \eqref{eq:abs2} with the outer approximation~\eqref{eq:outer1}-\eqref{eq:outer4}, one obtains an SOC-representable approximation of $H_\epsilon$. Note that this approximation is no longer conservative, but can be made so by instead using the polyhedral conservative approximation of $S_\epsilon$ as discussed in Section~\ref{sec:polyapprox}.

\subsection{Approximation via robust optimization}

An alternative conservative approximation which we consider is based on robust optimization~\cite{RobustBook}.
\begin{lemma}
Let $\epsilon \in (0,1)$ and suppose that $\xi$ follows an $n$-dimensional multivariate Guassian distribution. Without loss of generality, we assume each component is independent standard Gaussian with zero mean and unit variance. Let $\Gamma = F_{\chi_n}^{-1}(1-\epsilon)$, where $F_{\chi_n}^{-1}$ is the inverse cumulative distribution function of the chi distribution with $n$ degrees of freedom. Let
\begin{align}\label{eq:robustsdp}
R_\epsilon = \left\{(a,b,c,d,k,\lambda) : 
\left[\begin{array}{cccc}
\lambda I & & a & c \\
& k - \lambda\Gamma^2 & b & d\\
a^T & b& 1\\
c^T&d&&1
\end{array}\right] \succeq 0
 \right \},
\end{align}
where the notation $A \succeq 0$ means that the symmetric matrix $A$ is positive semidefinite, and blank entries represent zero blocks. Let $R_\epsilon^{proj}$ be the projection of the set $R_\epsilon$ onto the variables $(a,b,c,d,k)$. Then $R_\epsilon^{proj}$ is convex and $R_\epsilon^{proj} \subseteq H_\epsilon$. That is, the set $R_\epsilon^{proj}$ is a conservative, convex approximation of the chance constraint~\eqref{eq:quadprob}.
\begin{proof}
It is sufficient to show that if $(a,b,c,d,k,\lambda) \in R_\epsilon$, then there exists a set $U$ such that $P(\xi \in U) \ge 1-\epsilon$ such that
\begin{equation}\label{eq:robustquad}
(a^T\eta + b)^2 + (c^T\eta+d)^2 \leq k, \, \forall \, \eta \in U.
\end{equation}
Instead of allowing $U$ to vary for any point in the set, which is equivalent to the original chance constraint, we fix $U = \{\eta : ||\eta||_2 \le \Gamma\}$ and therefore obtain a conservative approximation. By the definition of the chi distribution, $P(\xi \in U) = 1-\epsilon$. In the terminology of robust optimization, $U$ is an uncertainty set. It is a standard result, which follows from the S-lemma and Schur complement lemmas, that $R_\epsilon^{proj}$ is precisely the set of points satisfying~\eqref{eq:robustquad} for this choice of the uncertainty set $U$~\cite{RobustBook}. Convexity follows since $R_\epsilon$ is the set of points satisfying a linear matrix inequality (LMI), which is tractable by semidefinite programming (SDP).
\end{proof}
\end{lemma}

Note that the choice of $\Gamma = F_{\chi_n}^{-1}(1-\epsilon)$ may be overly conservative, especially when $n$ is large, although we are not aware of any theoretical guidance on choosing a smaller value of $\Gamma$ such that the chance constraint remains satisfied.

\subsection{Nemirovski-Shapiro CVaR approximation}

The third approximation we consider is based on the so-called CVaR approximation proposed by Nemirovski and Shapiro~\cite{NS2007}. Let $I(z)$ be the indicator function of the interval $[0,\infty)$, i.e., $I(z) = 1$ if $z \ge 0$ and $I(z) = 0$ otherwise. We can rewrite the quadratic chance constraint in the following equivalent expected-value form:
\begin{equation}\label{eq:expvalform}
\mathbb{E}_\xi\left[I((a^T\xi + b)^2 + (c^T\xi+d)^2-k)\right] \le \epsilon.
\end{equation}

Nemirovski and Shapiro propose to upper bound the indicator function $I$ with the convex increasing function $\psi(z) = \max(1+z,0)$, which, up to rescaling ($z \to z/\alpha$ for some $\alpha$), is the best possible convex upper bound on the indicator function in the sense that if $\omega(z)$ is another convex increasing upper bound, then there exists $\alpha > 0$ such that $\psi(z/\alpha) \le \omega(z)$ for all $z \in \mathbb{R}$. Then the constraint
\begin{equation}
\inf_{\alpha > 0}\left[\mathbb{E}_\xi\left[\psi( ((a^T\xi + b)^2 + (c^T\xi+d)^2-k)/\alpha)\right] -\epsilon\right] \le 0,
\end{equation}
is a conservative approximation of the quadratic chance constraint~\eqref{eq:expvalform} which is furthermore convex in $(a,b,c,d,k)$, which motivates the following lemma.

\begin{lemma}
Let $\epsilon \in (0,1)$ and suppose that $\xi$ follows an $n$-dimensional multivariate Guassian distribution. Let
\begin{equation}
NS_\epsilon = \left\{(a,b,c,d,k,\alpha) : 
\mathbb{E}_\xi\left[\max((a^T\xi + b)^2 + (c^T\xi+d)^2-k+\alpha,0)\right] \le \alpha\epsilon, \alpha \ge 0 \right\}
\end{equation}
Let $NS_\epsilon^{proj}$ be the projection of the set $NS_\epsilon$ onto the variables $(a,b,c,d,k)$. Then $NS_\epsilon^{proj}$ is convex and $NS_\epsilon^{proj} \subseteq H_\epsilon$. That is, the set $NS_\epsilon^{proj}$ is a conservative, convex approximation of the chance constraint~\eqref{eq:quadprob}.
\begin{proof}
See~\cite{NS2007}.
\end{proof}
\end{lemma}

\subsection{A comparison of approximations}

We have presented three convex, conservative formulations of the quadratic chance constraint: one based on two-sided chance constraints, one based on robust optimization, and one based on convex approximation of the indicator function. All three have different tractability properties. In order of increasing computational difficulty, the two-sided approximation can be implemented, with small additional approximation error, by second-order cone programming (SOCP) following the developments presented in this work. The approximation based on robust optimization has an SDP formulation which may not be practical on large-scale problems, although we note the work of~\cite{ChordalSDP} where specialized methods were developed to exploit the block structure. The CVaR approximation is the most computationally challenging; it has no known reformulation in terms of standard problem classes and requires multidimensional integration to evaluate.

One might expect that the more computationally challenging approaches could yield tighter approximations. In this section, we examine a two-dimensional example in order to gain some understanding of the relative strengths of the approximations. We find, perhaps surprisingly, that no one approximation strictly dominates another. Hence, the two-sided approximation we propose has value in both its strength and ease of implementation.

As an example we will recall the simple case of
\begin{equation}\label{eq:nonconvex_example2}
\mathbb{P}((x\xi_1)^2 + (y\xi_2)^2 \le 1) \ge 1-\epsilon,
\end{equation}
where $\xi_1$ and $\xi_2$ are independent, standard Gaussian random variables.

Note that $||(\xi_1,\xi_2)||_2$ follows the chi distribution with 2 degrees of freedom, so in the robust approximation we can pick the uncertainty set $U = \{ (\eta_1,\eta_2) : ||(\eta_1,\eta_2)||_2 \leq F_{\chi_2}^{-1}(1-\epsilon) \}$ where $F_{\chi_2}^{-1}$ is the inverse cumulative distribution function of the chi distribution with two degrees of freedom.

In this example, \eqref{eq:robustsdp} reduces to
\begin{equation}\label{eq:psdexample}
\left[\begin{array}{ccccc}
\lambda &  & & x &  \\
  &\lambda &   &  & y\\
 &  & 1 - \lambda\Gamma^2 &   & \\
x & &  & 1 &  \\
&y&&&1
\end{array}\right] \succeq 0.
\end{equation}
By a Schur complement argument, the matrix~\eqref{eq:psdexample} is positive semidefinite iff $1-\lambda \Gamma^2 \ge 0, \lambda - x^2 \ge 0,$ and $\lambda - y^2 \ge 0$, which holds iff $x \in [-1/\Gamma,1/\Gamma]$ and $y \in [-1/\Gamma,1/\Gamma]$, a simple box constraint.

An interesting observation is that for $n = 2$, the choice of $\Gamma = F_{\chi_2}^{-1}(1-\epsilon)$ is \textit{minimal} in the sense that any smaller value no longer corresponds to a conservative approximation of the chance constraint~\eqref{eq:nonconvex_example2}:
\begin{equation}
\mathbb{P}(((1/\Gamma)\xi_1)^2 + ((1/\Gamma)\xi_2)^2 \le 1) = \mathbb{P}(\xi_1^2 + \xi_2^2 \le \Gamma^2) = F_{\chi_2}(\Gamma).
\end{equation}
In other words, the robust approximation to~\eqref{eq:nonconvex_example2} touches the boundary of the exact feasible set at the corners of the box. This observation eliminates the possibility of relaxing the overconservatism of the robust approximation by decreasing the size of the uncertainty set for the case of $n = 2$.

A point is feasible to the two-sided approximation \eqref{eq:quadapproxall} for $\beta = \frac{1}{2}$ iff $\exists f_1, f_2$ such that $f_1^2 + f_2^2 \leq 1$, $\mathbb{P}(|x\xi_1|\leq f_1) \geq 1-\frac{\epsilon}{2}$, and $\mathbb{P}(|y\xi_2|\leq f_2) \geq 1-\frac{\epsilon}{2}$.
By symmetry, these two chance constraints hold iff $f_1/|x| \geq \Phi^{-1}(1-\frac{\epsilon}{4})$ and $f_2/|y| \geq \Phi^{-1}(1-\frac{\epsilon}{4})$. Therefore, the point $(x,y)$ feasible to the two-sided approximation iff
\[
x^2 + y^2 \leq \frac{1}{\Phi^{-1}(1-\frac{\epsilon}{4})^2},
\]
a simple ball constraint.

The CVaR approximation, to our knowledge, does not yield a closed-form algebraic representation, although in this simple case we are able to evaluate it by numerical integration.

\begin{figure}[t]
\centering
\includegraphics[trim=75 10 75 10,clip,scale=0.7]{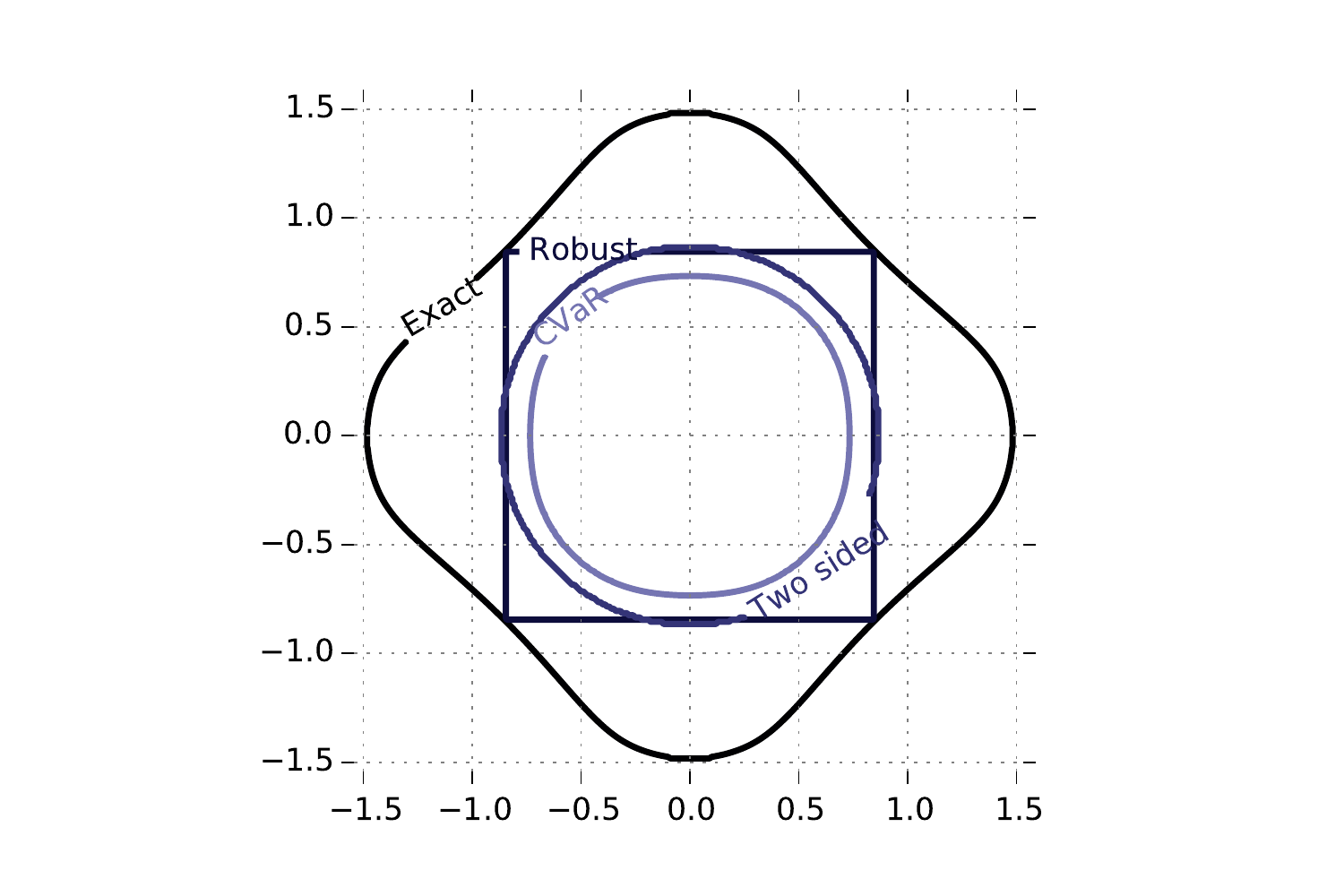}
\caption{ Outlined in black, the exact nonconvex feasible set $(x,y)$ satisfying $\mathbb{P}((x\xi_1)^2 + (y\xi_2)^2 \le 1) \ge 1-\epsilon$ for $\epsilon = 0.5$. We compare the three different convex approximations.}
\label{fig:eps05}
\end{figure}

\begin{figure}
\centering
\includegraphics[trim=75 10 75 10,clip,scale=0.6]{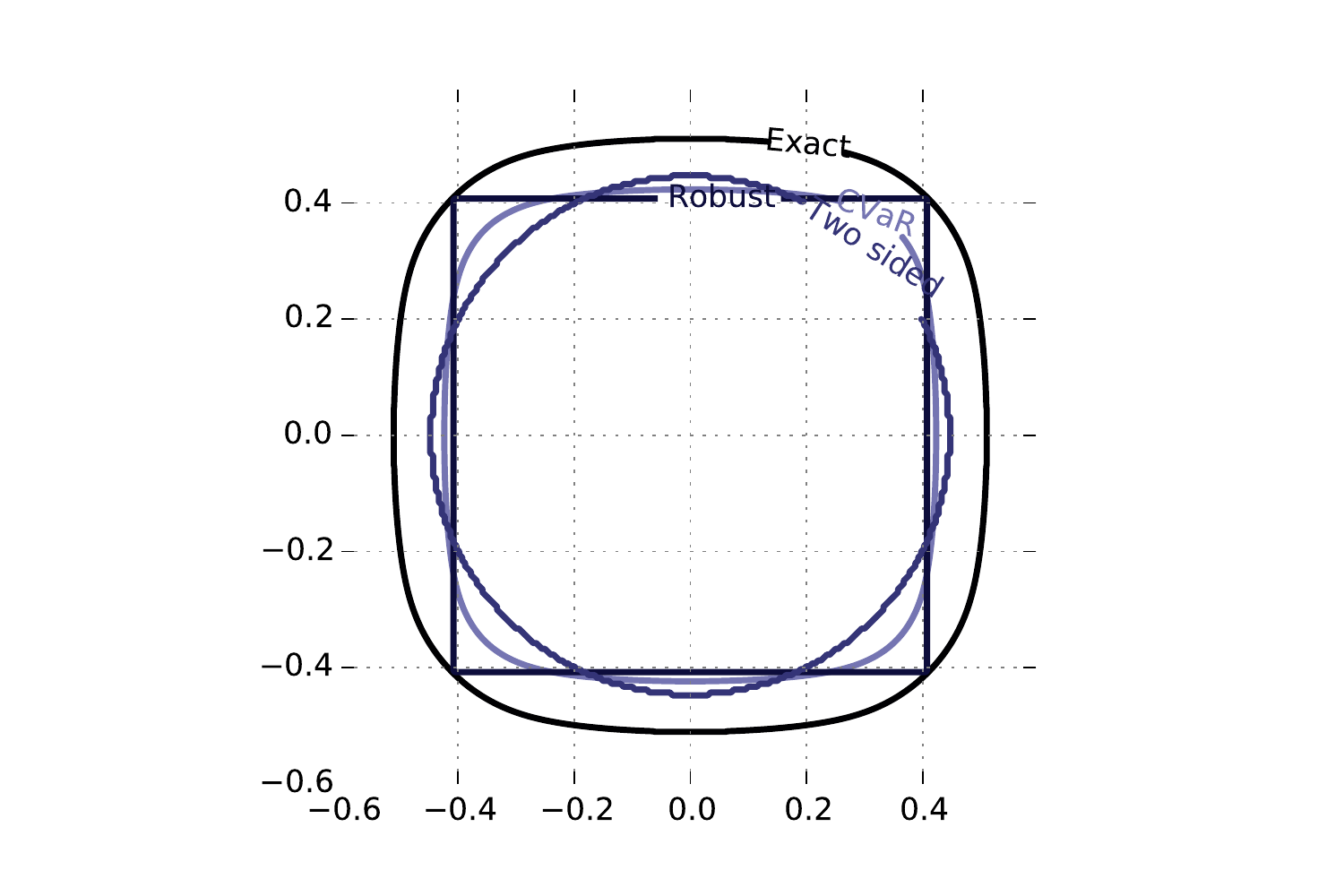}
\includegraphics[trim=75 10 75 10,clip,scale=0.6]{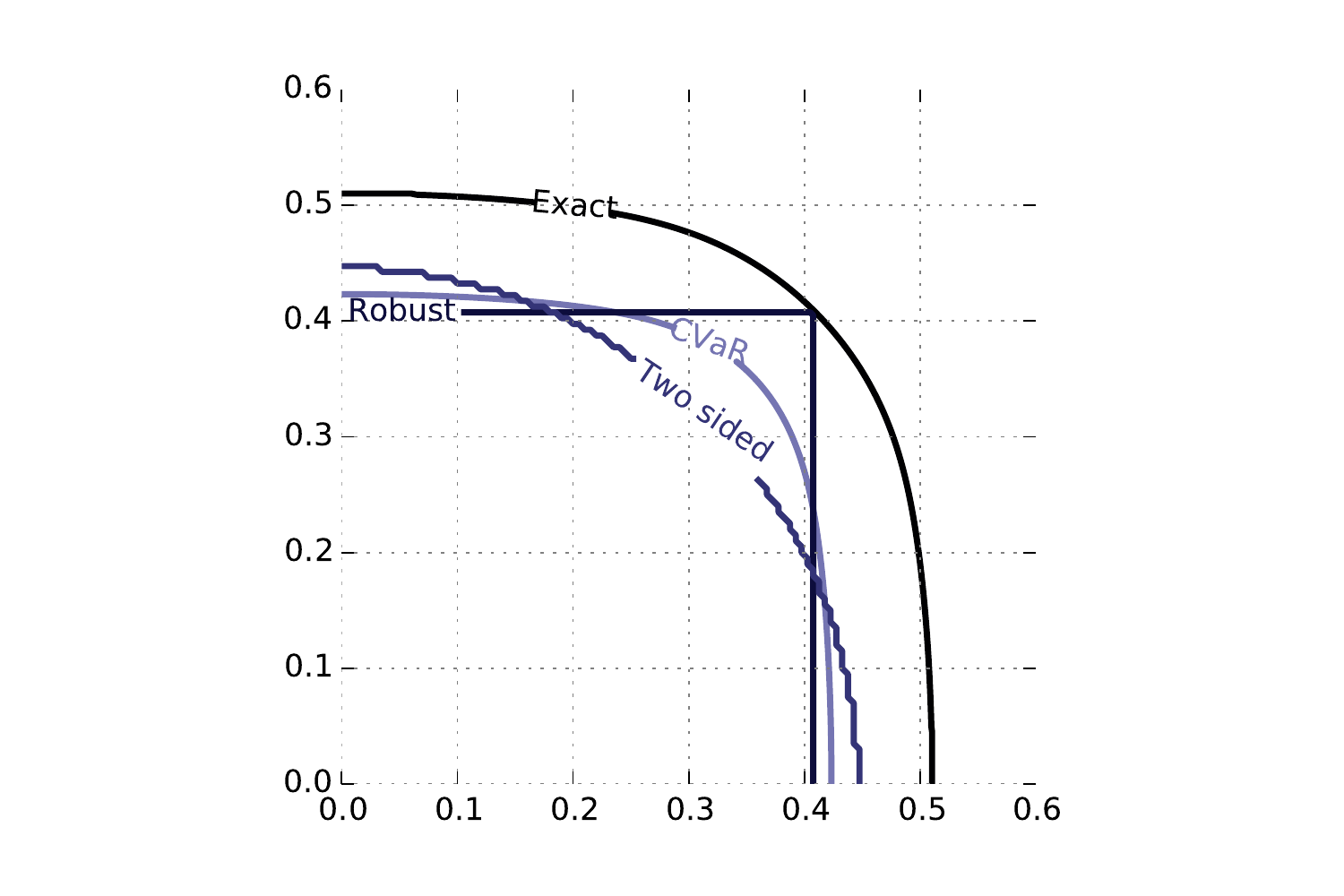}
\caption{ Outlined in black, the exact (seemingly convex) feasible set $(x,y)$ satisfying $\mathbb{P}((x\xi_1)^2 + (y\xi_2)^2 \le 1) \ge 1-\epsilon$ for $\epsilon = 0.05$. We compare the three different convex approximations. On the right, a zoomed-in view of the top-right corner shows that no approximation strictly dominates another. }
\label{fig:eps005}
\end{figure}

Figures~\ref{fig:eps05} and~\ref{fig:eps005} compare the three approximations with the exact feasible set for $\epsilon = 0.5$ and $\epsilon = 0.05$, respectively. For $\epsilon = 0.5$, both the robust and the two-sided approximations dominate the CVaR approximation. For $\epsilon = 0.05$, no approximation is a strict subset of another. Curiously, for this particular case the exact set $H_{0.05}$ appears to be convex.

\section{Conclusion}

Building on top of the basic convexity result for two-sided chance constraints developed in Section~\ref{sec:cvx2side}, we have shown, perhaps surprisingly, that a large class of more general nonlinear chance constraints is in fact convex (Theorem~\ref{thm:nonlinearchance}). In addition, our analysis of the computational tractability of the two-sided chance constraint, and in particular the polyhedral approximation of the set $S_\epsilon$ with provable approximation quality, develops practical methodologies which we believe are novel in the chance constraint literature. Finally, we have demonstrated that the two-sided chance constraint yields a useful approximation of the quadratic chance constraint which originally motivated this work.

We believe that our convexity results in Section~\ref{sec:cvx2side} can be easily extended to elliptical log-concave distributions following \cite{Lagoa05}. Extensions to more general distributions are not at all obvious, although the distributionally robust model in Section~\ref{sec:distrobust} may serve as a useful approximation. The conditions under which the quadratic chance constraint set $H_\epsilon$ is convex is left as an open question, although based on our computational experiments we conjecture that the set is convex for $\epsilon$ sufficiently small.

\section*{Acknowledgements}
We thank Michael (Misha) Chertkov of Los Alamos National Laboratory for discussions which inspired this work. M. Lubin was supported by the DOE Computational Science Graduate Fellowship,
which is provided under grant number DE-FG02-97ER25308.


\bibliography{refs}{}

\begin{thebibliography}{10}

\bibitem{ChordalSDP}
{\sc M.~Andersen, L.~Vandenberghe, and J.~Dahl}, {\em Linear matrix
  inequalities with chordal sparsity patterns and applications to robust
  quadratic optimization}, in Computer-Aided Control System Design (CACSD),
  2010 IEEE International Symposium on, Sept 2010, pp.~7--12.

\bibitem{RobustBook}
{\sc A.~Ben-Tal, L.~El~Ghaoui, and A.~Nemirovski}, {\em Robust Optimization},
  Princeton Series in Applied Mathematics, Princeton University Press, October
  2009.

\bibitem{BenTalNemirovskiRobust2000}
{\sc A.~Ben-Tal and A.~Nemirovski}, {\em Robust solutions of linear programming
  problems contaminated with uncertain data}, Mathematical Programming, 88
  (2000), pp.~411--424.

\bibitem{BergenBook}
{\sc A.~R. Bergen and V.~Vittal}, {\em Power systems analysis}, Prentice Hall,
  1999.

\bibitem{ccopf-sirev}
{\sc D.~Bienstock, M.~Chertkov, and S.~Harnett}, {\em Chance-constrained
  optimal power flow: Risk-aware network control under uncertainty}, SIAM
  Review, 56 (2014), pp.~461--495.

\bibitem{bolognani2015fast}
{\sc S.~Bolognani and F.~D{\"o}rfler}, {\em Fast power system analysis via
  implicit linearization of the power flow manifold}, in 53rd Annual Allerton
  Conference on Communication, Control, and Computing, 2015.

\bibitem{Bonmin}
{\sc P.~Bonami, L.~T. Biegler, A.~R. Conn, G.~Cornuéjols, I.~E. Grossmann,
  C.~D. Laird, J.~Lee, A.~Lodi, F.~Margot, N.~Sawaya, and A.~Wächter}, {\em An
  algorithmic framework for convex mixed integer nonlinear programs}, Discrete
  Optimization, 5 (2008), pp.~186 -- 204.

\bibitem{BoydBook}
{\sc S.~Boyd and L.~Vandenberghe}, {\em Convex Optimization}, Cambridge
  University Press, New York, NY, USA, 2004.

\bibitem{CharnesCooper}
{\sc A.~Charnes and W.~W. Cooper}, {\em Deterministic equivalents for
  optimizing and satisficing under chance constraints}, Operations Research, 11
  (1963), pp.~18--39.

\bibitem{Copulas}
{\sc J.~Cheng, M.~Houda, and A.~Lisser}, {\em Second-order cone programming
  approach for elliptically distributed joint probabilistic constraints with
  dependent rows}, 2014.
\newblock Available on Optimization Online.

\bibitem{working}
{\sc Y.~Dvorkin, L.~Roald, M.~Lubin, and M.~Chertkov}, {\em Chance constraints
  for improving the reliability of {ACOPF} solutions}, working paper, 2016.

\bibitem{HiriartLemarechal93book2}
{\sc J.~B. Hiriart-Urruty and C.~Lemar{\'e}chal}, {\em Convex Analysis and
  Minimization Algorithms}, vol.~I-II, Springer-Verlag, Germany, 1993.

\bibitem{OPFreview}
{\sc M.~Huneault and F.~Galiana}, {\em A survey of the optimal power flow
  literature}, Power Systems, IEEE Transactions on, 6 (1991), pp.~762--770.

\bibitem{Lagoa05}
{\sc C.~M. Lagoa, X.~Li, and M.~Sznaier}, {\em Probabilistically constrained
  linear programs and risk-adjusted controller design}, SIAM Journal on
  Optimization, 15 (2005), pp.~938--951.

\bibitem{JuMPChance}
{\sc M.~Lubin}, {\em {JuMPChance.jl}},  (2015).
\newblock {http://dx.doi.org/10.5281/zenodo.13740}.

\bibitem{JuMPChanceCaseStudy}
{\sc M.~Lubin, Y.~Dvorkin, and S.~Backhaus}, {\em A robust approach to chance
  constrained optimal power flow with renewable generation}, Power Systems,
  IEEE Transactions on, to appear (2015), pp.~1--10.

\bibitem{NS2007}
{\sc A.~Nemirovski and A.~Shapiro}, {\em Convex approximations of chance
  constrained programs}, SIAM Journal on Optimization, 17 (2007), pp.~969--996.

\bibitem{PrekopaBook}
{\sc A.~Pr\'{e}kopa}, {\em Stochastic Programming}, Springer Netherlands, 1995.

\bibitem{schrijver2003}
{\sc A.~Schrijver}, {\em Combinatorial Optimization: Polyhedra and Efficiency},
  no.~v. 1 in Algorithms and Combinatorics, Springer, 2003.

\bibitem{StubbsMehrotra99}
{\sc R.~A. Stubbs and S.~Mehrotra}, {\em A branch-and-cut method for 0-1 mixed
  convex programming}, Mathematical Programming, 86 (1999), pp.~515--532.

\bibitem{VanAckooij10}
{\sc W.~Van~Ackooij, R.~Henrion, A.~M\"{o}ller, and R.~Zorgati}, {\em On
  probabilistic constraints induced by rectangular sets and multivariate normal
  distributions}, Mathematical Methods of Operations Research, 71 (2010),
  pp.~535--549.

\end{thebibliography}

\bibliographystyle{siam}

\end{document}